\documentclass[11pt,a4paper,twoside]{article}

\newcommand{\theTitle}{Exact solutions and bounds for network SIR and SEIR models using a rooted-tree approximation}
\newcommand{\theShortTitle}{Rooted-tree approximation for network contagion}  
\newcommand{\theAuthors}{C. L. Hall and B. A. Siebert}
\newcommand{\theShortAuthors}{CL Hall and BA Siebert}

\title{\theTitle}
\author{\theAuthors}
\date{\today}

\usepackage[a4paper,top=2cm,bottom=2cm,left=2.5cm,right=2.5cm,marginparwidth=2cm,includeheadfoot,twoside]{geometry}
\usepackage{fancyhdr}
\usepackage{setspace}
\usepackage{parskip}

\usepackage{amsmath}
\usepackage{amsfonts}
\usepackage{amssymb}

\usepackage{graphicx}
\graphicspath{{./}{./Figures/}}
\usepackage{caption}
\usepackage{subcaption}

\usepackage{tikz}
\usepackage{pgfplots}
\pgfplotsset{compat=1.15}
\usetikzlibrary{calc}
\usetikzlibrary{plotmarks}
\usetikzlibrary{patterns}
\usetikzlibrary{shapes}
\usetikzlibrary{external}
\tikzsetexternalprefix{./Figures/}

\usepackage{color}

\usepackage[english]{babel}

\usepackage{lipsum}               % Filler text
\usepackage{hhline}               % Better lines in tables
\usepackage{todonotes}

\makeatletter
\renewcommand{\todo}[2][]{\tikzexternaldisable\@todo[#1]{#2}\tikzexternalenable}
\makeatother

\newcommand{\diff}[2]{\frac{\mathrm{d} #1}{\mathrm{d} #2}}

\newcommand{\ie}{\emph{i.e.}}
\newcommand{\eg}{\emph{e.g.}}

\definecolor{myGrey}{RGB}{240,240,240}

\def\XXint#1#2#3{{\setbox0=\hbox{$#1{#2#3}{\int}$}
     \vcenter{\hbox{$#2#3$}}\kern-.5\wd0}}

\newcommand{\Sus}[1]{\langle S_{#1} \rangle}
\newcommand{\SusApprox}[1]{\langle S^{*}_{#1} \rangle}

\newcommand{\Exp}[2]{\langle E_{#1}^{(#2)} \rangle}
\newcommand{\ExpVec}[1]{\langle \mathbf{E}_{#1} \rangle}
\newcommand{\ExpVecApprox}[1]{\langle \mathbf{E}^{*}_{#1} \rangle}
\newcommand{\ExpSimp}[1]{\langle E_{#1} \rangle}
\newcommand{\ExpVecInit}{\langle \mathbf{E}_0 \rangle^\text{init}}
\newcommand{\ExpSimpInit}{\langle E_0 \rangle^\text{init}}
\newcommand{\Inf}[1]{\langle I_{#1} \rangle}
\newcommand{\InfApprox}[1]{\langle I^{*}_{#1} \rangle}
\newcommand{\InfInit}{\langle I_0 \rangle^\text{init}}
\newcommand{\Rec}[1]{\langle R_{#1} \rangle}
\newcommand{\RecApprox}[1]{\langle R^{*}_{#1} \rangle}

\newcommand{\QuuVec}[1]{\langle \mathbf{Q}_{#1} \rangle}

\newcommand{\QuuVecApprox}[1]{\langle \mathbf{Q}^{*}_{#1} \rangle}

\newcommand{\InfSus}[2]{\langle I_{#1} S_{#2} \rangle}
\newcommand{\SusSus}[2]{\langle S_{#1} S_{#2} \rangle}
\newcommand{\ExpSus}[3]{\langle E_{#1}^{(#2)} S_{#3} \rangle}

\newcommand{\RecSus}[2]{\langle R_{#1} S_{#2} \rangle}
\newcommand{\InfRecSus}[3]{\langle I_{#1} R_{#2} S_{#3} \rangle}

\newcommand{\classu}{u}
\newcommand{\classv}{v}
\newcommand{\nodek}{k}
\newcommand{\nodej}{j}
\newcommand{\nodei}{i}
\newcommand{\iindex}{n}
\newcommand{\vecx}{\mathbf{x}}
\newcommand{\vecxk}[1]{x_{#1}}
\newcommand{\vecf}{\mathbf{f}}
\newcommand{\vecfk}[1]{f_{#1}}

\newcommand{\phijn}[2]{{\varphi}^{(#1)}_{#2}}
\newcommand{\phijnp}[3]{{\varphi}^{(#1)}_{#2 \leftarrow #3}}
\newcommand{\phivecn}[1]{\boldsymbol{\varphi}_{#1}}
\newcommand{\phin}[1]{{\varphi}_{#1}}
\newcommand{\lamtotn}[1]{\lambda_{#1}}
\newcommand{\lamtotnp}[2]{\lambda_{#1 \leftarrow #2}}
\newcommand{\ajkn}[3]{b^{(#1 \leftarrow #2)}_{#3}}
\newcommand{\mujn}[2]{\mu^{(#1)}_{#2}}
\newcommand{\muvecn}[1]{\boldsymbol{\mu}_{#1}}
\newcommand{\musimpn}[1]{{\mu}_{#1}}
\newcommand{\nujn}[2]{\nu^{(#1)}_{#2}}
\newcommand{\nuvecn}[1]{\boldsymbol{\nu}_{#1}}
\newcommand{\nusimpn}[1]{{\nu}_{#1}}
\newcommand{\gamman}[1]{\gamma_{#1}}
\newcommand{\An}[1]{\mathbf{B}_{#1}}
\newcommand{\Mn}[1]{\mathbf{M}_{#1}}

\newcommand{\Cn}[1]{C_{#1}}
\newcommand{\netState}{\mathbf{X}}
\newcommand{\nodeState}[1]{X_{#1}}

\newcommand{\numClasses}{N_u}
\newcommand{\numNodes}{N}
\newcommand{\zerovec}{\boldsymbol{0}}
\newcommand{\uvec}{\mathbf{e}}
\newcommand{\parn}[1]{p(#1)}
\newcommand{\pparn}[1]{p[p(#1)]}
\newcommand{\Neigh}{\mathcal{N}}

\usepackage{hyperref}
\hypersetup{
	unicode=true,                 % non-Latin characters in Acrobat bookmarks
	pdftoolbar=true,              % show Acrobat toolbar?
	pdfmenubar=true,              % show Acrobat menu?
	pdffitwindow=true,            % page fit to window when opened
	pdftitle={\theTitle},         % title
	pdfauthor={\theShortAuthors}, % author
	pdfsubject={},                % subject of the document
	pdfnewwindow=true,            % links in new window
	pdfkeywords={},               % list of keywords
	colorlinks=true,              % false: boxed links; true: colored links
	linkcolor=black,              % color of internal links
	citecolor=black,              % color of links to bibliography
	filecolor=black,              % color of file links
	urlcolor=black,               % color of external links
}

\usepackage[noabbrev,capitalise]{cleveref}

\begin{document}
	
\maketitle

\begin{abstract}
	
	In this paper, we develop a node-based approximate model for Markovian contagion dynamics on networks. 
	We prove that our approximate model is exact for SIR (susceptible-infectious-recovered) and SEIR (susceptible-exposed-infectious-recovered) dynamics on tree graphs with a single source of infection and that the model otherwise gives upper bounds on the probabilities of each node being susceptible.
	Our analysis of SEIR contagion dynamics is generalised to SEIR models with arbitrarily many distinct classes of exposed state.
	In the case of trees with a single source of infection, our approach yields a system of partially-decoupled linear differential equations that exactly describes the evolution of node-state probabilities. 
	We use this to state explicit closed-form solutions for SIR dynamics on a chain.
\end{abstract}

\section{Introduction}
\label{S:Introduction}

Network-based models have been used extensively to describe the spread of a contagious state through a population via the connections between individuals. 
Such models are particularly important in describing the spread of disease \cite{Danon2011,KissMillerSimonEpidemics,Miller2014,PastorSatorras2015} but have also been used to study social contagion \cite{Hill2010,Zadeh2014,Zanette2002}, financial contagion \cite{Gai2010,Leventides2019}, and cascading failure in power systems \cite{Guo2017,Zhang2017}.
Network contagion is also a rich field for theoreticians;
exact solutions to network contagion models are mostly unavailable and so both efficient numerical methods and good approximate models are valuable \cite{PastorSatorras2015}.

In many contagion models on networks, each node represents an individual and each edge represents a contact or connection that facilitates the spread of contagion between nodes.
At any given time, each node has a state (\textit{e.g.}, susceptible, infectious,%
%\footnote{We use the term `infectious' rather than the standard `infected' throughout to improve the clarity of our analysis of SEIR models.
%In this analysis, it is convenient for us to use the phrase `become infected' to describe both S to E and S to I transitions; as a result, our terminology is clearer if we refer to the I state as `infectious' rather than `infected'.}
or recovered in the classic SIR model \cite{NewmanNetworks,PastorSatorras2015}) and the node states evolve over time according to the rules that constitute the contagion model.
%We refer to the state of all nodes in a network at a specified time as the \textit{network state}; for a contagion model with $\numNodeStates$ types of node state and a network with $\numNodes$ nodes, there are $\numNodeStates^\numNodes$ different network states.
In many such models, node state evolution is probabilistic and occurs over continuous time; 
in these cases, the spread of contagion through the network is a continuous-time discrete-space stochastic process where the state space is the set of states for all nodes in the network.

%Let $\{\netState(t)\}$ represent the stochastic process for network contagion dynamics on a network of $\numNodes$ nodes.
%Any realisation of this process can be represented as a time-dependent $\numNodes$-dimensional vector of node states, $\netState(t)$, so that $\nodeState{\nodek}(t)$ gives the state of the $\nodek$th node at time $t$. 
%%We use roman letters to indicate these node states; for an SIR model, we have $X_\nodek(t) \in \{\text{S},\text{I},\text{R}\}$.
%Following various other authors \cite{Sharkey2015,Sharkey2015a,Simon2018,Wilkinson2014}, we use angle brackets to indicate probabilities. 
%We define $\Sus{\nodek}(t) = P\left[X_k(t) = \text{S}\right]$ to be the probability that node $\nodek$ is susceptible at time $t$, and we define $\Inf{\nodek}$ and $\Rec{\nodek}$, similarly.
%We also define $\InfSus{\nodej}{\nodek}(t) = P\left[X_j(t) = \text{I} \cap X_k(t) = \text{S}\right]$ to be the probability that node $\nodej$ is infectious and node $\nodek$ is susceptible at time $t$, and we define other joint probabilities similarly. 

One challenge with stochastic network contagion models is to determine the node state probabilities as functions of time. %for a given network with a given initial network state.
%the probability that a given node is in a given state at a given time after a given initial network state.
Even for very simple contagion models, this is difficult on large networks because the node states do not evolve independently.
In the most general case, node state probabilities can only be determined exactly from network state probabilities, which in turn can only be determined exactly by solving the master equations for the stochastic process. 
Since the size of the state space increases geometrically with the number of nodes,% (and therefore so does the number of master equations), 
this is not computationally feasible on any but the smallest networks.

Instead, various methods have been developed for estimating---and, in some cases, bounding---node state probabilities in network contagion models.
The simplest of these is the node-based mean field approximation \cite{PastorSatorras2015} (also called the first-order model \cite{NewmanNetworks}, the individual-based model \cite{KissMillerSimonEpidemics}, or the $N$-intertwined mean field approximation \cite{Simon2018,VanMieghem2009}). 
In this approach, node state probabilities are assumed to be independent of each other, so that joint probabilities can be expressed as the product of individual node state probabilities.
%;
%using angle-bracket notation for probabilities (see \cite{Sharkey2015,Sharkey2015a,Simon2018,Wilkinson2014}) this can be written as $\InfSus{\nodej}{\nodek} = \Inf{\nodej}\Sus{\nodek}$.
While this is a useful assumption that closes the evolution equations for node state probabilities, it is not perfectly accurate. 
In reality, the states of neighbouring nodes are positively correlated: \eg, the neighbours of a susceptible node are more likely to be susceptible than would be expected from assuming independence \cite{Donnelly1993}.
As a result, the node-based mean field approximation applied to standard contagion models will typically overestimate rates of infection and hence underestimate the probability that a given node is susceptible. %$\Sus{\nodek}$ \cite{Cator2018}.

Two other approaches used to estimate and bound node state probabilities are the pair-based approximation \cite{Cator2012,PastorSatorras2015} and the message passing approximation \cite{Karrer2010}.
To develop the pair-based approximation, Cator and Van Mieghem \cite{Cator2012} introduced variables for the joint probabilities of the states of neighbouring nodes and they derived evolution equations for these probabilities using a closure approximation to exclude the dependence on higher-order moments.
To develop the message-passing approximation, Karrer and Newman \cite{Karrer2010} considered the directed edges of the network and developed expressions for the probabilities that infection has not yet been transmitted along each edge. 

While these two approaches are conceptually very different, Wilkinson and Sharkey \cite{Wilkinson2014} showed that they are equivalent for Markovian SIR dynamics. 
Pair-based and message-passing approximations are more computationally demanding than node-based approximations but are generally more accurate than the node-based mean field approximation.
When the underlying network is a tree, both approaches yield exact results for the SIR model \cite{Karrer2010,Sharkey2015}.

In this paper, we develop and analyse a new approximate model of network contagion that can be applied to Markovian SIR and SEIR (susceptible-exposed-infectious-recovered) contagion models, including SEIR models with multiple distinct exposed states. 
The approximation we derive is a `node-based' approximation; it takes the form of a closed system of differential equations for node state probabilities.
As such, our approximation has a similar level of computational complexity to the node-based mean field approximation and is considerably simpler than the pair-based or message-passing approximations.

We refer to our approximation as the `rooted-tree approximation' because it yields exact results on trees with a single initially-infectious node.
This contrasts with both the node-based mean field model, which can never give exact results, and the pair-based and message-passing approximations, which give exact results on any tree regardless of the number of initially-infectious nodes \cite{Karrer2010,Sharkey2015}.
The exact differential equations obtained using our approximation are very simple and lead to explicit closed-form solutions for node state probabilities on rooted trees. 
We believe that these explicit solutions have not previously been reported.

On other networks (non-trees or trees with multiple initially-infectious nodes), we prove that the rooted-tree approximation gives \textit{upper} bounds on the probabilities that nodes are susceptible.
This contrasts with the other approximations described above, which give \textit{lower} bounds on the probabilities that nodes are susceptible;
this lower bound result is generally understood to hold for node-based mean field approximation of SIR models \cite{Cator2018,Cator2014,Donnelly1993} and
has been proved for node-based mean field approximation of SIS models \cite{Donnelly1993,Simon2018}
and for pair-based/message-passing approximation of SIR models \cite{Karrer2010,Wilkinson2014}

%the node-based mean field approximation, which has been proven to give lower bounds on the probabilities that nodes are susceptible in an SIS model \cite{Donnelly1993,Simon2018} and is generally understood to give lower bounds on these probabilities in SIR models \cite{Cator2018,Cator2014,Donnelly1993}.
%Similarly, our result contrasts with the pair-based and message-passing approximations, which have been proven to give lower bounds on the probabilities that nodes are susceptible in SIR models \cite{Karrer2010,Wilkinson2014}.
%Since the node-based mean field approximation will typically also give a lower bound on these probabilities, our work suggests that the two node-based approximations  could be combined to yield new, more accurate node-based models. 

The development of our approximation exploits the fact that neither the SIR nor SEIR models permit the possibility of reinfection. 
In the case of an SIR model on a tree with a single initially-infectious node, this enables us to formulate an exact expression for the rate of infection in terms of the probabilities that nodes are susceptible.
For other networks and initial conditions, a similar approach enables us to formulate a cooperative system of differential equations %(\ie, a system that satisfies the Kamke--M{\"u}ller conditions \cite{Donnelly1993,Simon2018}) 
where the approximate rate of infection is a lower bound on the true rate of infection.
This enables us to use Simon and Kiss's methods from \cite{Simon2018} to prove that our approach yields upper bounds on the probablilities that nodes are susceptible.

%The key results of this paper are the rooted-tree approximations given in system XX\todo{EQUATION REFS} for SIR dynamics and YY for generalised SEIR dynamics. 

%Our paper is structured as follows. 
%In \cref{S:SIR} we present our development of the rooted-tree approximation for the network SIR model.
%We begin by showing how our approach leads to an exact solution on a rooted tree (\cref{S:RootedTree-SIR}) that yields simple closed-form solutions to node state probabilities (\cref{S:Solutions-SIR}).
%We then further develop our analysis by showing that the rooted-tree approximation also gives an upper bound for the probability that a node is susceptible in a general graph (\cref{S:Bounds-SIR}).
%In \cref{S:SEIR} we then repeat this analysis for a generalised SEIR model with arbitrarily many exposed states; 
%we find that an analogous approach to \cref{S:RootedTree-SIR} yields further exact solutions on rooted trees (\cref{S:RootedTree-SEIR})  and we again show that the rooted-tree approximation yields a bound on the true node state probabilities (\cref{S:Bounds-SEIR}).

Our main contribution in this paper can be summarised as the rooted-tree approximation systems given in \eqref{E:RootedTreeApprox-SIR} and \eqref{E:RootedTreeApprox-SEIR} for SIR and SEIR models respectively.
In \cref{S:SIR}, we develop \eqref{E:RootedTreeApprox-SIR} for SIR models and prove that it is exact on rooted trees and otherwise yields an upper bound on the probability of being susceptible.
In \cref{S:SEIR}, we repeat this analysis for SEIR models to develop \eqref{E:RootedTreeApprox-SEIR}.
Finally, in \cref{S:Discussion}, we discuss the merits and limitations of our approach and make comparisons with other theoretical approaches to network contagion.
We conclude by offering avenues for further exploration and extension of the rooted-tree approximation.

%Our paper is structured as follows. 
%In \cref{S:SIR} we present our development of the rooted-tree approximation for the network SIR model.
%We begin by showing how our approach leads to an exact solution on a rooted tree (\cref{S:RootedTree-SIR}) that yields simple closed-form solutions to node state probabilities (\cref{S:Solutions-SIR}).
%We then further develop our analysis by showing that the rooted-tree approximation also gives an upper bound for the probability that a node is susceptible in a general graph (\cref{S:Bounds-SIR}).
%In \cref{S:SEIR} we then repeat this analysis for a generalised SEIR model with arbitrarily many exposed states; 
%we find that an analogous approach to \cref{S:RootedTree-SIR} yields further exact solutions on rooted trees (\cref{S:RootedTree-SEIR})  and we again show that the rooted-tree approximation yields a bound on the true node state probabilities (\cref{S:Bounds-SEIR}).
%
%Our main contributions can be summarised as the rooted-tree approximations given in \eqref{E:RootedTreeApprox-SIR} and \eqref{E:RootedTreeApprox-SEIR} for SIR and SEIR models respectively.
%In \cref{S:Discussion}, we discuss the merits and limitations of our approach and how it compares to other theoretical approaches to network contagion.
%We conclude by offering avenues for further exploration and extension of the rooted-tree approximation.

%NOTE FOR CONCLUSION OF PAPER 2: Discuss extension to semi-Markov models? (rate depends on time in disease state?)

\section{Rooted-tree approximation for the SIR model}
\label{S:SIR}

\subsection{Preliminaries}
\label{S:Prelims-SIR}

Let $\{\netState(t)\}$ represent the stochastic process for network contagion dynamics on a network of $\numNodes$ nodes.
Any realisation of this process can be represented as a time-dependent $\numNodes$-dimensional vector of node states, $\netState(t)$, so that $\nodeState{\nodek}(t)$ gives the state of the $\nodek$th node at time $t$. 
%We use roman letters to indicate these node states; for an SIR model, we have $X_\nodek(t) \in \{\text{S},\text{I},\text{R}\}$.
Following various other authors \cite{Sharkey2015,Sharkey2015a,Simon2018,Wilkinson2014}, we use angle brackets to indicate probabilities. 
Specifically, we define $\Sus{\nodek}(t) = P\left[X_k(t) = \text{S}\right]$ to be the probability that node $\nodek$ is susceptible at time $t$, 
we define $\InfSus{\nodej}{\nodek}(t) = P\left[X_j(t) = \text{I} \cap X_k(t) = \text{S}\right]$ to be the probability that node $\nodej$ is infectious and node $\nodek$ is susceptible at time $t$, and we define other probabilities and joint probabilities similarly.
%and we define $\Inf{\nodek}$ and $\Rec{\nodek}$, similarly.
%We also define $\InfSus{\nodej}{\nodek}(t) = P\left[X_j(t) = \text{I} \cap X_k(t) = \text{S}\right]$ to be the probability that node $\nodej$ is infectious and node $\nodek$ is susceptible at time $t$, and we define other joint probabilities similarly.

In this section, we focus on the standard network SIR model as described in \cite{NewmanNetworks} and elsewhere. 
At any time, each node can either be susceptible (S), infectious (I) or recovered (R) and node states change over time according to a Markovian process.
Susceptible nodes in contact with infectious nodes become infected at rate $\lamtotn{}$; 
that is, the probability that a susceptible node in contact with an infectious node becomes infectious in the next $\Delta t$ is given by $\lamtotn{} \Delta t + o(\Delta t)$.
Infection rates are taken to be additive over neighbours, so that additional infectious neighbours will increase the probability that a susceptible node becomes infectious in a given $\Delta t$.
Infectious nodes recover at rate $\gamman{}$ regardless of the states of their neighbours.

As a further generalisation, we assume that $\lamtotn{}$ can depend on the associated directed edge, and that $\gamman{}$ can depend on the associated node.
Thus, we assume that the rate of infection can depend on the nodes involved %, that this is not symmetric (\ie, the rate at which node $\nodej$ would infect node $\nodek$ can be different from the rate at which node $\nodek$ would infect node $\nodej$), 
and that the rate of recovery from infection can vary from node to node.
We represent this using subscripts, so that $\lamtotnp{\nodek}{\nodej}$ is the rate at which node $\nodek$ becomes infected given that node $\nodek$ is susceptible and node $\nodej$ is infectious, and $\gamman{\nodek}$ is the rate at which node $\nodek$ would recover given that it is currently infectious.

With this notation, the following is an exact description of node probability dynamics for an SIR model on a network:
\begin{subequations}
\label{E:GeneralEvol-SIR}
\begin{align}
	\diff{\Sus{\nodek}}{t}    &= - \sum_{j \in \Neigh(k)} \lamtotnp{\nodek}{\nodej} \InfSus{\nodej}{\nodek}, \label{E:SusEvol-SIR}\\
	\diff{\Inf{\nodek}}{t}    &= \sum_{j \in \Neigh(k)} \lamtotnp{\nodek}{\nodej} \InfSus{\nodej}{\nodek} - \gamman{\nodek} \Inf{\nodek},  \label{E:InfEvol-SIR}\\
	\diff{\Rec{\nodek}}{t}    &= \gamman{\nodek} \Inf{\nodek}, \label{E:RecEvol-SIR}
\end{align}
\end{subequations}
where $\Neigh(\nodek)$ represents the set of upstream neighbours of node $\nodek$ (\ie, the set of nodes $\nodej$ for which $\lamtotnp{\nodek}{\nodej}$ is nonzero).

%This system of equations is not closed; in order to construct a node-based model of contagion dynamics, we need to supplement this system with an expression for the pair probability $\InfSus{\nodej}{\nodek}$ in terms of node state probabilities. 
%In \cref{S:RootedTree-SIR}, we describe how this can be done exactly for a tree with a single initially-infectious node.

\subsection{Exact SIR dynamics on a rooted tree}
\label{S:RootedTree-SIR}

Consider the case where the underlying network is a tree and where a single node is infectious at $t = 0$ and all other nodes are susceptible.
We assign the the label $\nodek = 0$ to the initially-infectious node and identify it as the root of the tree.
We will use the term `rooted tree' throughout our analysis (including for SEIR models) to refer to a tree where there is a unique node that is not in a susceptible or recovered state at $t = 0$.
For any other node $\nodek \neq 0$, it is possible to identify a unique parent node $\parn{\nodek}$ as the neighbour of $\nodek$ that lies between node $\nodek$ and the root.
Since all infection spreads from the root node % and the only path to node $\nodek$ from the root node is via node $\parn{\nodek}$, 
it follows that node $\nodek$ can only be infected by node $\parn{\nodek}$.
This enables us to simplify our notation and analysis in this section: we define $\lamtotn{\nodek} = \lamtotnp{\nodek}{\parn{\nodek}}$ as the rate at which node $\nodek$ is infected by its parent node, and we omit the sums in equations \eqref{E:SusEvol-SIR} and \eqref{E:InfEvol-SIR}.

Thus, the evolution of node state probabilities on a rooted tree is given by
\begin{subequations}
\label{E:RootedEvol-SIR}
\begin{align}
	\diff{\Sus{\nodek}}{t}    &= 
	\begin{cases}
			0, & \nodek = 0, \\
			-\lamtotn{\nodek} \InfSus{\parn{\nodek}}{\nodek}, & \nodek \neq 0;
	\end{cases}
	\label{E:SusEvol-SIR-Rooted}\\
	\diff{\Inf{\nodek}}{t}    &= 
	\begin{cases}
		- \gamman{\nodek} \Inf{\nodek}, & \nodek = 0, \\
		\lamtotn{\nodek} \InfSus{\parn{\nodek}}{\nodek} - \gamman{\nodek} \Inf{\nodek}, & \nodek \neq 0;
	\end{cases}
	\label{E:InfEvol-SIR-Rooted}\\
	\diff{\Rec{\nodek}}{t}    &= \gamman{\nodek} \Inf{\nodek}. \label{E:RecEvol-SIR-Rooted}
\end{align}	
\end{subequations}
These equations need to be solved subject to initial conditions
\begin{align}
	\Sus{\nodek}(0) &= 
	\begin{cases}
		0, & \nodek=0, \\
		1, & \nodek\neq 0;
	\end{cases}, &
	\Inf{\nodek}(0) &= 
	\begin{cases}
	1, & \nodek=0, \\
	0, & \nodek\neq 0;
	\end{cases}, &	
	\Rec{\nodek}(0) &= 0. \label{E:RootedICs-SIR}
\end{align}
This system of equations is not closed; in order to construct a node-based model of contagion dynamics, we need expressions for the pair probabilities $\InfSus{\parn{\nodek}}{\nodek}$ in terms of the node state probabilities. 
The analysis below shows how this can be achieved exactly.

Consider any node $\nodek \neq 0$.
The law of total probability gives
\begin{equation}
	\Sus{\nodek}
	= \SusSus{\parn{\nodek}}{\nodek} 
	+ \InfSus{\parn{\nodek}}{\nodek} 
	+ \RecSus{\parn{\nodek}}{\nodek}
	\label{E:SusExhaustionSIR}
\end{equation}
Since infection can only spread from node $\parn{\nodek}$ to node $\nodek$ and not \textit{vice versa}, we find that $\nodeState{\parn{\nodek}} = \text{S}$ implies $\nodeState{\nodek} = \text{S}$ (\ie, if the parent of node $k$ is susceptible then node $k$ must also be susceptible).
Hence, $\SusSus{\parn{\nodek}}{\nodek} = \Sus{\parn{\nodek}}$
and \eqref{E:SusExhaustionSIR} can be rearranged as
\begin{equation}
	\InfSus{\parn{\nodek}}{\nodek} = \Sus{\nodek} -  \Sus{\parn{\nodek}} - \RecSus{\parn{\nodek}}{\nodek}.
	\label{E:InfSusFromRecSus-SIR}
\end{equation}
This indicates that an expression for $\RecSus{\parn{\nodek}}{\nodek}$ in terms of node state probabilities could be used to obtain an expression for $\InfSus{\parn{\nodek}}{\nodek}$ in terms of node state probabilities.

We note that the only way to achieve a state where $X_{\parn{\nodek}} = \text{R}$ and $X_{\nodek} = \text{S}$ is for node $\parn{\nodek}$ to recover  while node $\nodek$ is susceptible. 
Once such a state is achieved, it will persist permanently since node $\parn{\nodek}$ will remain recovered and node $\nodek$ cannot become infected except via node $\parn{\nodek}$.
Expressed mathematically, this means that
\begin{equation}
	\diff{\RecSus{\parn{\nodek}}{\nodek}}{t} =  \gamman{\parn{\nodek}} \InfSus{\parn{\nodek}}{\nodek}, \qquad \nodek \neq 0,
\end{equation}
which can be rearranged using \eqref{E:SusEvol-SIR-Rooted} to yield
\begin{equation}
	\diff{\RecSus{\parn{\nodek}}{\nodek}}{t}    = - \frac{\gamman{\parn{\nodek}}}{\lamtotn{\nodek}} \diff{\Sus{\nodek}}{t}, \qquad  \nodek \neq 0. \label{E:RecSusEvol-SIR}
\end{equation}
Integrating \eqref{E:RecSusEvol-SIR} and applying the initial conditions $\RecSus{\parn{\nodek}}{\nodek}(0) = 0$ and $\Sus{\nodek}(0) = 1$ for $\nodek \neq 0$, we find that
%\begin{equation}
$	\RecSus{\parn{\nodek}}{\nodek} = \frac{\gamman{\parn{\nodek}}}{\lamtotn{\nodek}} - \frac{\gamman{\parn{\nodek}}}{\lamtotn{\nodek}} \Sus{\nodek}$.
	%\label{E:RecSus-SIR}
%\end{equation}
Substituting this into \eqref{E:InfSusFromRecSus-SIR} then yields
\begin{equation}
	\InfSus{\parn{\nodek}}{\nodek} = \frac{\lamtotn{\nodek}+\gamman{\parn{\nodek}}}{\lamtotn{\nodek}} \Sus{\nodek} -  \Sus{\parn{\nodek}} - \frac{\gamman{\parn{\nodek}}}{\lamtotn{\nodek}}, \qquad \nodek \neq 0.
	\label{E:InfSus-SIR}
\end{equation}

Equation \eqref{E:InfSus-SIR} gives an expression for $\InfSus{\parn{\nodek}}{\nodek}$ purely in terms of the node state probabilities $\Sus{\nodek}$ and  $\Sus{\parn{\nodek}}$.
Substituting into system \eqref{E:RootedEvol-SIR}, we obtain the following closed system for the node state probabilities:
\begin{subequations}
\label{E:RootedEvol-SIR-Full}
\begin{align}
	\diff{\Sus{\nodek}}{t}    &= 
	\begin{cases}
		0, & \nodek = 0, \\
		-(\lamtotn{\nodek}+\gamman{\parn{\nodek}}) \Sus{\nodek} + \lamtotn{\nodek} \Sus{\parn{\nodek}} + \gamman{\parn{\nodek}}, & \nodek \neq 0;
	\end{cases}
	\label{E:SusEvol-SIR-RootedFull}\\
	\diff{\Inf{\nodek}}{t}    &= 
	\begin{cases}
		- \gamman{\nodek} \Inf{\nodek}, & \nodek = 0, \\
		(\lamtotn{\nodek}+\gamman{\parn{\nodek}}) \Sus{\nodek} - \lamtotn{\nodek} \Sus{\parn{\nodek}} - \gamman{\parn{\nodek}} - \gamman{\nodek} \Inf{\nodek}, & \nodek \neq 0;
	\end{cases}
	\label{E:InfEvol-SIR-RootedFull}\\
	\diff{\Rec{\nodek}}{t}    &= \gamman{\nodek} \Inf{\nodek}. \label{E:RecEvol-SIR-RootedFull}
\end{align}
\end{subequations}
This system can be solved subject to the initial conditions in \eqref{E:RootedICs-SIR} to yield an exact representation of node state probabilities on a rooted tree.

\subsection{Closed form solutions}
\label{S:Solutions-SIR}

The system in \eqref{E:RootedEvol-SIR-Full} is amenable to further analysis leading to explicit closed form solutions.
We observe that the differential equations in \eqref{E:RootedEvol-SIR-Full} are all linear and have constant coefficients.
Moreover, the system is partially decoupled:  
the equations for $\diff{\Sus{\nodek}}{t}$ are independent of $\Inf{\nodek}$ and $\Rec{\nodek}$, 
the equations for $\diff{\Inf{\nodek}}{t}$ are independent of $\Rec{\nodek}$, 
and all equations for node state probabilities at a given node are independent of the states of the node's children and siblings.
It follows that the differential equations in \eqref{E:RootedEvol-SIR-Full} can be solved sequentially using standard methods for first-order constant coefficients linear differential equations.

For example, consider the case where $\lamtotn{}$ and $\gamman{}$ are constant for all nodes.
In this case, the symmetry of the system implies that node state probabilities will be identical for nodes of equal depth (\ie, equal distance from the root node).
Thus, we can obtain all node state probabilities by considering a chain of nodes labelled 0, 1, 2, \textit{etc.} where each node is connected to its ordinal neighbours.

Rearranging \eqref{E:RootedEvol-SIR-Full} and exploiting the fact that $\Sus{\nodek} + \Inf{\nodek} + \Rec{\nodek} = 1$, the system to be solved for this `chain' problem is
\begin{subequations}
\label{E:ChainEvol-SIR}
\begin{align}
		\diff{\Sus{\nodek}}{t} + (\lamtotn{}+\gamman{}) \Sus{\nodek}  &= 
	\begin{cases}
		0, & \nodek = 0, \\
		\lamtotn{} \Sus{\nodek-1} + \gamman{}, & \nodek \neq 0;
	\end{cases}
	\label{E:SusEvol-SIR-Chain}\\
	\diff{\Inf{\nodek}}{t}  + \gamman{} \Inf{\nodek}  &= 
	-\diff{\Sus{\nodek}}{t};
	\label{E:InfEvol-SIR-Chain}\\
	\Rec{\nodek}   &= 1 - \Sus{\nodek} - \Inf{\nodek}, \label{E:RecEvol-SIR-Chain}
\end{align}
\end{subequations}
subject to the initial conditions \eqref{E:RootedICs-SIR}.

This system can be solved explicitly using a range of different methods (\eg, operator $D$ methods or Laplace transforms). %in \cref{S:SolutionMethod} we use operator $D$ methods to show that
Applying any of these solution methods, we find that
%using operator $D$ methods or Laplace transforms, we find that
\begin{subequations}
	\label{E:ClosedFormSolution}
\begin{align}
	\Sus{\nodek}(t) &= 1 
		- \frac{\lamtotn{}^\nodek}{(\lamtotn{} + \gamman{})^\nodek} 
		+ \frac{\lamtotn{}^\nodek}{(\lamtotn{} + \gamman{})^\nodek} \mathrm{e}^{-(\lamtotn{} + \gamman{})t}
		\sum_{\iindex=0}^{\nodek-1} \frac{(\lamtotn{} + \gamman{})^\iindex t^\iindex}{\iindex!}, 
		\label{E:ClosedFormSolution-S} \\
	\Inf{\nodek}(t) &= \mathrm{e}^{-\gamman{}t} 
		- \mathrm{e}^{-(\lamtotn{} + \gamman{})t} \sum_{\iindex=0}^{\nodek-1} \frac{\lamtotn{}^\iindex t^\iindex}{\iindex!}, 
		\label{E:ClosedFormSolution-I} \\
	\Rec{\nodek}(t) &= \frac{\lamtotn{}^\nodek}{(\lamtotn{} + \gamman{})^\nodek}
	- \mathrm{e}^{-\gamman{}t}
	+ \mathrm{e}^{-(\lamtotn{} + \gamman{})t} \sum_{\iindex=0}^{\nodek-1} \left[\left(\lamtotn{}^{\iindex} - \frac{\lamtotn{}^\nodek}{(\lamtotn{}+\gamman{})^{\nodek-\iindex}}\right)  \frac{t^{\iindex}}{\iindex!} \right].
	\label{E:ClosedFormSolution-R}
\end{align}	
\end{subequations}
To the best of our knowledge, this is the first time that this simple, closed-form solution has been reported in the literature on contagion on networks.

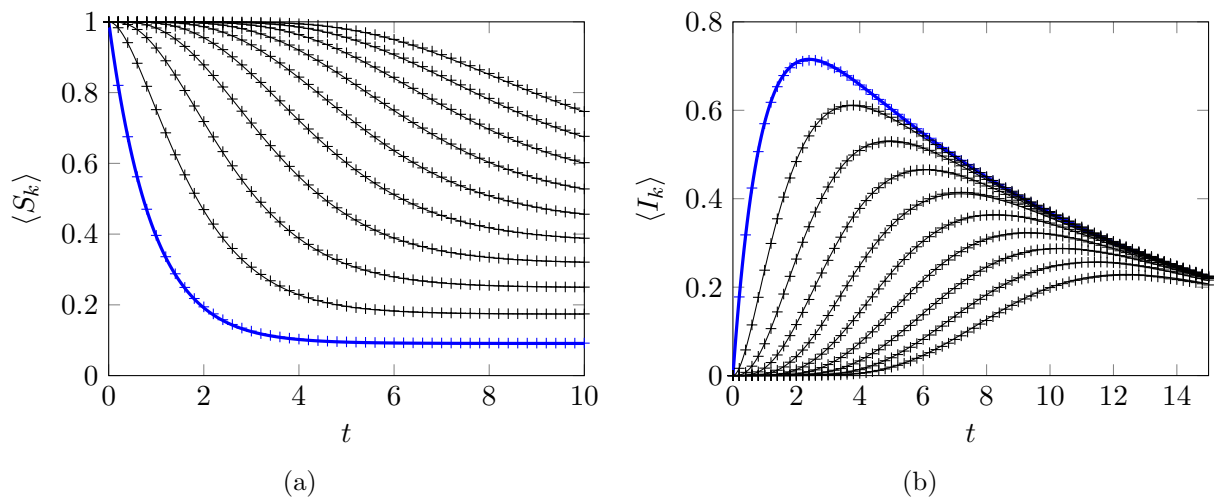
\begin{figure}
	\centering
	\begin{subfigure}{0.49\textwidth}
		\centering
		\tikzsetnextfilename{SIR-Chain-SusAllNodes}
		\begin{tikzpicture}
			\pgfplotstableread{fig-SIR-Chain-SusAllNodes-DE-A.txt}\deSolutionTable
			\pgfplotstableread{fig-SIR-Chain-SusAllNodes-Gil-A.txt}\gilSolutionTable
			\begin{axis}[
				ymin=0,ymax=1,xmin=0,xmax=10,
				xlabel={$t$},
				ylabel={$\Sus{\nodek}$},
				%	xtick={0},
				%	ytick={0},
				%	extra x ticks={0.5,1},
				%	extra y ticks={0.5,1},
				%	extra x tick labels={$\tfrac{1}{2}$, $1$}, 
				%	extra y tick labels={$\tfrac{1}{2}$, $1$}, 
				width={\textwidth},
				height={0.8\textwidth}
				]
				\addplot[color=blue, very thick] table[x index = 0, y index = 1] from \deSolutionTable;
				\addplot[color=blue, mark=+, only marks] table[x index = 0, y index = 1] from \gilSolutionTable;
				\foreach \figLineMaker in {2,3,...,10} {
					\addplot[color=black] table[x index = 0, y index = \figLineMaker] from \deSolutionTable;
					\addplot[color=black, mark=+, only marks] table[x index = 0, y index = \figLineMaker] from \gilSolutionTable;
				}
			\end{axis}
		\end{tikzpicture}
%		\caption{Probability of nodes being susceptible at time $t$.}
		\caption{}
		\label{fig:SIR-Chain-SusAllNodes}
	\end{subfigure}%
	\hspace{\fill}
	\begin{subfigure}{0.49\textwidth}
		\centering
		\tikzsetnextfilename{SIR-Chain-InfAllNodes}
		\begin{tikzpicture}
			\pgfplotstableread{fig-SIR-Chain-InfAllNodes-DE-A.txt}\deSolutionTable
			\pgfplotstableread{fig-SIR-Chain-InfAllNodes-Gil-A.txt}\gilSolutionTable
			\begin{axis}[
				ymin=0,ymax=0.8,xmin=0,xmax=15,
				xlabel={$t$},
				ylabel={$\Inf{\nodek}$},
				%	xtick={0},
				%	ytick={0},
				%	extra x ticks={0.5,1},
				%	extra y ticks={0.5,1},
				%	extra x tick labels={$\tfrac{1}{2}$, $1$}, 
				%	extra y tick labels={$\tfrac{1}{2}$, $1$}, 
				width={\textwidth},
				height={0.8\textwidth}
				]
				\addplot[color=blue, very thick] table[x index = 0, y index = 1] from \deSolutionTable;
				\addplot[color=blue, mark=+, only marks] table[x index = 0, y index = 1] from \gilSolutionTable;
				\foreach \figLineMaker in {2,3,...,10} {
					\addplot[color=black] table[x index = 0, y index = \figLineMaker] from \deSolutionTable;
					\addplot[color=black, mark=+, only marks] table[x index = 0, y index = \figLineMaker] from \gilSolutionTable;
				}
			\end{axis}
		\end{tikzpicture}
%		\caption{Probability of nodes being infectious at time $t$.} 
		\caption{}
\label{fig:SIR-Chain-InfAllNodes}
		%\label{TK}
	\end{subfigure}
	\caption{\small 
		Comparision of the rooted-tree solutions for $\Sus{\nodek}$ and $\Inf{\nodek}$ in \eqref{E:ClosedFormSolution} with simulation results from the average of $10^5$ Gillespie algorithm simulations of the full stochastic SIR model. 
		Subfigure (a) shows results for $\Sus{\nodek}$ while subfigure (b) shows results for $\Inf{\nodek}$.
		In both cases, the rooted tree solutions are shown as continuous lines and the numerical results are shown as points marked $+$.
		Results are shown for the first ten nodes; results from $\nodek = 1$ are indicated with a thicker blue line and subsequent nodes produce curves further to the right.
		Parameters used are $\lamtotn{} = 1$ and $\gamman{} = 0.1$.} \label{fig:SIR-Chain-SusInfAllNodes}
\end{figure}

\cref{fig:SIR-Chain-SusInfAllNodes} shows comparisons of $\Sus{\nodek}(t)$ and $\Inf{\nodek}(t)$ from \eqref{E:ClosedFormSolution} with empirical node state probabilities based on averaging $10^5$ Gillespie algorithm simulations of the underlying stochastic model.
All calculations were performed in \textsc{Matlab} and code is provided at \url{https://github.com/cameronlhall/rootedtreeapprox}.
These figures illustrate the fact that \eqref{E:ClosedFormSolution} are exact results; the theoretical results for  $\Sus{\nodek}(t)$ and $\Inf{\nodek}(t)$ are virtually indistinguishable from results obtained using Gillespie simulations.

\cref{fig:SIR-Chain-SusInfAllNodes} also illustrates some properties of SIR dynamics on a chain that can be derived from analysis of \eqref{E:ClosedFormSolution}.
For example, \eqref{E:ClosedFormSolution-I} can be rearranged as
\begin{equation}
	\Inf{\nodek}(t) = \mathrm{e}^{-\gamman{}t}\left(1  
	- \mathrm{e}^{-\lamtotn{}t} \sum_{\iindex=0}^{\nodek-1} \frac{\lamtotn{}^\iindex t^\iindex}{\iindex!} \right). \label{E:ClosedFormSolution-IAlt}
\end{equation}
Since the sum in \eqref{E:ClosedFormSolution-IAlt} is the first $k$ terms in the Maclaurin series of $\mathrm{e}^{\lambda t}$, we see that $\Inf{k}(t)$ will initially be close to zero and will remain close to zero for longer for larger values of $\nodek$.
Additionally, we observe that the term in brackets in \eqref{E:ClosedFormSolution-IAlt} will asymptotically approach $1$ as $t \to \infty$, which implies that $\Inf{\nodek}(t) \sim \mathrm{e}^{-\gamman{}t}$ as $t \to \infty$.
Both the early time behaviour where $\Inf{\nodek}$ is close to zero and the late time behaviour where $\Inf{\nodek}\sim \mathrm{e}^{-\gamman{}t}$ are visible in \cref{fig:SIR-Chain-InfAllNodes}.

While \eqref{E:ChainEvol-SIR} and \eqref{E:ClosedFormSolution} are simple and elegant results, they are of limited practical use because they are specific to rooted trees.
Results that only hold on trees are not useful for describing contagion on contact networks or social networks because such networks tend to be highly clustered \cite{NewmanNetworks} and the clustering coefficient of a tree is necessarily zero.
However, \eqref{E:ChainEvol-SIR} can be adapted to obtain a node-based approximation of contagion dynamics that gives a bound on $\Sus{\nodek}$ for all networks.

\subsection{Bounds for SIR dynamics on a general network}
\label{S:Bounds-SIR}

In \cref{S:RootedTree-SIR}, we showed that the closed system \eqref{E:RootedEvol-SIR-Full} is equivalent to the system \eqref{E:RootedEvol-SIR}, which describes the evolution of node state probabilities for SIR dynamics on a rooted tree.
In this section, we develop an analogue of \eqref{E:RootedEvol-SIR-Full} that can be applied to a general network.
We show that this new formulation yields upper bounds on the functions $\Sus{\nodek}(t)$.

We begin by assuming that no node is recovered at $t = 0$, and so we can specify initial conditions where $\Sus{\nodek}(0)$ is given for each node and
\begin{align}
	\Inf{\nodek}(0) &= 1 - \Sus{\nodek}(0), & \Rec{\nodek}(0) &= 0.
\end{align}
We make this assumption without loss of generality since the recovered state is permanent in the SIR model; SIR dynamics on a network with initially-recovered nodes will be equivalent to SIR dynamics on a network where those nodes and associated edges have been removed.

%Our analysis will focus on system \eqref{E:GeneralEvol-SIR}, replicated below:
%\begin{align*}
%	\diff{\Sus{\nodek}}{t}    &= - \sum_{j \in \Neigh(k)} \lamtotnp{\nodek}{\nodej} \InfSus{\nodej}{\nodek}, 
%	%\label{E:SusEvol-SIR}
%	\\
%	\diff{\Inf{\nodek}}{t}    &= \sum_{j \in \Neigh(k)} \lamtotnp{\nodek}{\nodej} \InfSus{\nodej}{\nodek} - \gamman{\nodek} \Inf{\nodek},  
%	%\label{E:InfEvol-SIR}
%	\\
%	\diff{\Rec{\nodek}}{t}    &= \gamman{\nodek} \Inf{\nodek}, 
%	%\label{E:RecEvol-SIR}
%\end{align*}

%Our analysis will be applicable for all initial conditions.

%%Since any node that i
%
%Since the recovered state is permanent, SIR dynamics on a network where some nodes are recovered at $t = 0$ will be equivalent to SIR dynamics on an identical network where the recovered nodes and associated edges have been removed.
%%This means that a system in which some nodes have a nonzero probability of being recovered at $t=0$ can be analysed by taking a weighted average of the results from systems where initially-recovered nodes have been removed.

The analysis that follows is analogous to the derivation of the exact solution for rooted trees in \cref{S:RootedTree-SIR}, but we derive inequalities throughout.
Let $\nodej$ and $\nodek$ be chosen so that $\nodej \in \Neigh(\nodek)$.
From the laws of probability, we note that $\SusSus{\nodej}{\nodek} + \InfSus{\nodej}{\nodek} + \RecSus{\nodej}{\nodek} = \Sus{\nodek}$, and that $\SusSus{\nodej}{\nodek} \leq \Sus{\nodej}$.
Combining these gives
\begin{equation}
	\InfSus{\nodej}{\nodek} \geq \Sus{\nodek} - \Sus{\nodej} - \RecSus{\nodej}{\nodek}. \label{E:InfSusInequal-SIR}
\end{equation}

Now consider the dynamics of $\RecSus{\nodej}{\nodek}$. 
We note that a state where $\nodeState{\nodej} = \text{R}$ and $\nodeState{\nodek} = \text{S}$ can only arise from a state where node $\nodeState{\nodej} = \text{I}$ and $\nodeState{\nodek} = \text{S}$.
Additionally, a state where $\nodeState{\nodej} = \text{R}$ and $\nodeState{\nodek} = \text{S}$ can change to another state only if node $\nodek$ becomes infected from one of its neighbours. Thus,
\begin{equation}
	\diff{\RecSus{\nodej}{\nodek}}{t} = \gamman{\nodej} \InfSus{\nodej}{\nodek} - \sum_{\nodei \in \Neigh(\nodek)} \lamtotnp{\nodek}{\nodei} \InfRecSus{\nodei}{\nodej}{\nodek},
\end{equation}
and, since all probabilities are nonnegative, it follows that
\begin{equation}
	\diff{\RecSus{\nodej}{\nodek}}{t} \leq \gamman{\nodej} \InfSus{\nodej}{\nodek}. \label{E:RecSusInequal-SIR}
\end{equation}

Noting that the terms inside the summation in \eqref{E:SusEvol-SIR} are all nonnegative, we observe that
%\begin{equation}
%	-\diff{\Sus{\nodek}}{t} = \sum_{\nodej \in \Neigh(\nodek)} \lamtotnp{\nodek}{\nodej} \InfSus{\nodej}{\nodek},
%\end{equation}
%and hence
\begin{equation}
	-\diff{\Sus{\nodek}}{t} \geq \lamtotnp{\nodek}{\nodej} \InfSus{\nodej}{\nodek}. \label{E:SusEvolInequal-SIR}
\end{equation}
Combining \eqref{E:RecSusInequal-SIR} and \eqref{E:SusEvolInequal-SIR} then gives
\begin{equation}
	\diff{\RecSus{\nodej}{\nodek}}{t} \leq - \frac{\gamman{\nodej}}{\lamtotnp{\nodek}{\nodej}} \diff{\Sus{k}}{t}.
	\label{E:RecSusDE-Temp}
\end{equation}
Using the assumption that no nodes are recovered at $t = 0$, we recall that  $\RecSus{j}{k}(0) = 0$. 
This enables us to integrate \eqref{E:RecSusDE-Temp} from $t = 0$ to obtain
$
%\begin{equation}
	\RecSus{\nodej}{\nodek}(t) \leq \frac{\gamman{\nodej}}{\lamtotnp{\nodek}{\nodej}} \left[\Sus{\nodek}(0) - \Sus{\nodek}(t)\right],
$
%\end{equation}
and hence \eqref{E:InfSusInequal-SIR} becomes
\begin{equation}
	\InfSus{\nodej}{\nodek} \geq \Sus{\nodek} - \Sus{\nodej} - \frac{\gamman{\nodej}}{\lamtotnp{\nodek}{\nodej}} \left[\Sus{\nodek}(0) - \Sus{\nodek}(t)\right].
\end{equation}

Since $\InfSus{\nodej}{\nodek}$ is also nonnegative, it follows that
\begin{equation}
	\InfSus{\nodej}{\nodek}(t) \geq \left[\Sus{\nodek}(t) - \Sus{\nodej}(t) - \frac{\gamman{\nodej}}{\lamtotnp{\nodek}{\nodej}} \left[\Sus{\nodek}(0) - \Sus{\nodek}(t)\right] \right]^{+}, \label{E:MainInfSusInequal-SIR}
\end{equation}
where $[x]^{+}$ is defined so that
\begin{equation}
	[x]^{+} =\begin{cases}
		0, & x \leq 0, \\
		x, & x > 0.
	\end{cases} \label{E:PosPartDefn-SIR}
\end{equation}
Substituting into \eqref{E:SusEvol-SIR}, we obtain
\begin{equation}
	\diff{\Sus{\nodek}}{t} \leq  - \sum_{\nodej \in \Neigh(\nodek)} 
	\Big[
	- \gamman{\nodej} \Sus{\nodek}(0)
	+ (\lamtotnp{\nodek}{\nodej} + \gamman{\nodej}) \Sus{\nodek}(t)
	- \lamtotnp{\nodek}{\nodej} \Sus{\nodej}(t)
	\Big]^{+}. \label{E:MainSusEvolInequal-SIR}
\end{equation}

The differential inequality \eqref{E:MainSusEvolInequal-SIR} holds for the true node state probabilities $\Sus{\nodek}(t)$.
Based on this inequality, we now consider the relationship between the true solutions $\Sus{\nodek}(t)$ and approximate solutions $\SusApprox{\nodek}(t)$ that satisfy the system
\begin{equation}
	\diff{\SusApprox{\nodek}}{t} =  - \sum_{\nodej \in \Neigh(\nodek)} 
	\Big[
	- \gamman{\nodej} \SusApprox{\nodek}(0)
	+ (\lamtotnp{\nodek}{\nodej} + \gamman{\nodej}) \SusApprox{\nodek}(t)
	- \lamtotnp{\nodek}{\nodej} \SusApprox{\nodej}(t)
	\Big]^{+}, \label{E:ApproxSusEvol-SIR}
\end{equation}
subject to initial conditions
\begin{equation}
	\SusApprox{\nodek}(0) = \Sus{\nodek}(0). \label{E:ApproxSusICs-SIR}
\end{equation}

We will show that $\SusApprox{\nodek}(t) \geq \Sus{\nodek}(t)$ for all $\nodek$ and for all $t$.
This follows from the application of Lemma 1 from Simon and Kiss \cite{Simon2018}.
In order to use this result, we need to show that \eqref{E:ApproxSusEvol-SIR} is a cooperative system of differential equations.
This can be done using the Kamke--M{\"{u}}ller sufficient conditions \cite{Donnelly1993,Simon2018}, which state that an autonomous system
$
	\diff{\vecx}{t} = \vecf(\vecx),
$ 
will be cooperative as long as $\vecfk{\nodek}$ is a nondecreasing function of $\vecxk{\nodej}$ for all $\nodej \neq \nodek$.
In our case, we define $\vecx$ so that $\vecxk{\nodek} = \SusApprox{\nodek}$, and we define $\vecf(\vecx)$ so that
\begin{equation}
	\vecfk{\nodek}(\vecx) = 
	- \sum_{\nodej \in \Neigh(\nodek)} 
	\Big[
	- \gamman{\nodej} \SusApprox{\nodek}(0)
	+ (\lamtotnp{\nodek}{\nodej} + \gamman{\nodej}) \vecxk{\nodek}
	- \lamtotnp{\nodek}{\nodej} \vecxk{\nodej}
	\Big]^{+}.
\end{equation}
Since $\vecfk{\nodek}(\vecx)$ is continuous and the constants $\lamtotnp{\nodek}{\nodej}$ are nonnegative, it is clear that $\vecfk{\nodek}$ is a nondecreasing function of $\vecxk{\nodej}$ for all $\nodej$. 
Hence, the Kamke--M{\"{u}}ller conditions are satisfied and \eqref{E:ApproxSusEvol-SIR} is a cooperative system.
Using this fact alongside the initial conditions in \eqref{E:ApproxSusICs-SIR}, we apply Lemma 1 from \cite{Simon2018} to conclude that $\SusApprox{\nodek}(t) \geq \Sus{\nodek}(t)$ for all $\nodek$ and for all $t$.

To summarise this result, we can combine \eqref{E:ApproxSusEvol-SIR} with an equation for $\InfApprox{\nodek}$ based on \eqref{E:InfEvol-SIR} to obtain
\begin{subequations}
	\label{E:RootedTreeApprox-SIR}
\begin{align}
		\diff{\SusApprox{\nodek}}{t} &=  - \sum_{\nodej \in \Neigh(\nodek)} 
	\Big[
	- \gamman{\nodej} \SusApprox{\nodek}(0)
	+ (\lamtotnp{\nodek}{\nodej} + \gamman{\nodej}) \SusApprox{\nodek}(t)
	- \lamtotnp{\nodek}{\nodej} \SusApprox{\nodej}(t)
	\Big]^{+}, \label{E:RootedTreeApprox-SIR-Sus}\\
	\diff{\InfApprox{\nodek}}{t} &= \sum_{\nodej \in \Neigh(\nodek)} 
	\Big[
	- \gamman{\nodej} \SusApprox{\nodek}(0)
	+ (\lamtotnp{\nodek}{\nodej} + \gamman{\nodej}) \SusApprox{\nodek}(t)
	- \lamtotnp{\nodek}{\nodej} \SusApprox{\nodej}(t)
	\Big]^{+} -\gamman{\nodek}\InfApprox{\nodek}(t). \label{E:RootedTreeApprox-SIR-Exp}
\end{align}
\end{subequations}
If we also introduce $\RecApprox{\nodek} = 1 - \SusApprox{\nodek} - \InfApprox{\nodek}$, this gives a closed system of equations for the approximate dynamics of all node state probabilities.
We refer to system \eqref{E:RootedTreeApprox-SIR} as the rooted-tree approximation for SIR dynamics.

If the underlying network is a rooted tree, we can show that \eqref{E:RootedTreeApprox-SIR} is equivalent to \eqref{E:RootedEvol-SIR-Full}.
To see this, we note that $\SusApprox{\nodek}(t) \leq \SusApprox{\nodek}(0)$ for all time and that  $\SusApprox{\nodej}(t) \geq \SusApprox{\nodek}(t)$ for any $\nodej \in \Neigh(\nodek)$ other than $\nodej = \parn{\nodek}$.
Hence, the terms inside the square brackets in \eqref{E:RootedTreeApprox-SIR} will be nonpositive for any $\nodej \neq \parn{\nodek}$ and applying the positive part operator yields \eqref{E:RootedEvol-SIR-Full}.
If a network is known to be a rooted tree but the root is not identified, \eqref{E:RootedTreeApprox-SIR} will yield an exact solution without it being necessary to compute the parent of each node, as would be needed in order to use \eqref{E:RootedEvol-SIR-Full}.
%It follows that the rooted-tree approximation in \eqref{E:RootedTreeApprox-SIR} gives an upper bound on $\Sus{\nodek}(t)$ that is exact for the case of a rooted tree.

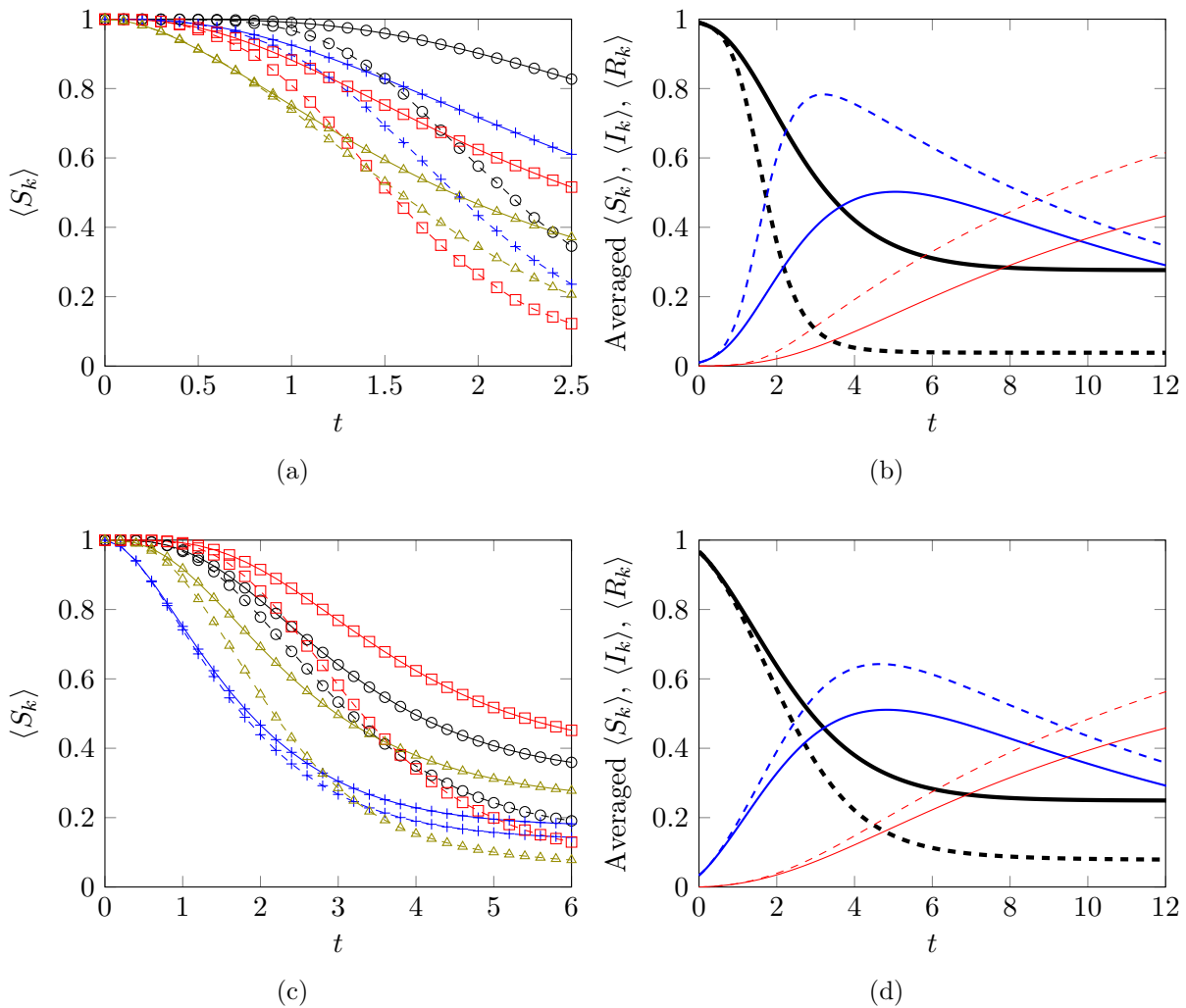
\begin{figure}
	\centering
	\begin{subfigure}[b]{0.49\textwidth}
	\centering
	\tikzsetnextfilename{SIR-ER-SusSelectedNodes}
	\begin{tikzpicture}
		\pgfplotstableread{fig-SIR-ER-SusSelectedNodes-All-D.txt}\mainTable
		\begin{axis}[
	ymin=0,ymax=1,xmin=0,xmax=2.5,
	xlabel={$t$},
	ylabel={$\Sus{\nodek}$},
%		xtick={0,0.05,0.1,0.15,0.2},
%		xticklabels={0,0.05,0.1,0.15,0.2},
%		xtick={0,0.05,0.1,0.15,0.2,0.25,0.3},
%		xticklabels={0,0.05,0.1,0.15,0.2,0.25,0.3},
	%	ytick={0},
	%	extra x ticks={0.5,1},
	%	extra y ticks={0.5,1},
	%	extra x tick labels={$\tfrac{1}{2}$, $1$}, 
	%	extra y tick labels={$\tfrac{1}{2}$, $1$}, 
	width={\textwidth},
	height={0.8\textwidth}
	]
	\addplot[color=black,mark=o] table[x index = 0, y index = 1] from \mainTable;
	\addplot[color=black,mark=o,dashed,mark options=solid] table[x index = 0, y index = 6] from \mainTable;
	\addplot[color=blue,mark=+] table[x index = 0, y index = 2] from \mainTable;
	\addplot[color=blue,mark=+,dashed,mark options=solid] table[x index = 0, y index = 7] from \mainTable;
	\addplot[color=red,mark=square] table[x index = 0, y index = 3] from \mainTable;
	\addplot[color=red,mark=square,dashed,mark options=solid] table[x index = 0, y index = 8] from \mainTable;
	\addplot[color=olive,mark=triangle] table[x index = 0, y index = 4] from \mainTable;
	\addplot[color=olive,mark=triangle,dashed,mark options=solid] table[x index = 0, y index = 9] from \mainTable;
%	\foreach \figLineMaker in {2,3,...,6} {
%		\addplot[color=black] table[x index = 0, y index = \figLineMaker] from \mainTable;
%	}
%	\foreach \figLineMaker in {7,8,...,12} {
%		\addplot[color=black, dashed] table[x index = 0, y index = \figLineMaker] from \mainTable;
%	}
\end{axis}
	\end{tikzpicture}
%	\caption{
%		Probabilities of selected nodes (indicated using different colours and mark styles) in an Erd\H{o}s--R\'{e}nyi graph being susceptible at time $t$, according to the rooted-tree approximation (continuous lines) and the average of $10^5$ Gillespie simulations (dashed lines).} 
\caption{}
	\label{fig:SIR-ER-SusSelectedNodes}
\end{subfigure}
\begin{subfigure}[b]{0.49\textwidth}
	\centering
	\tikzsetnextfilename{SIR-ER-AvgSusInf}
	\begin{tikzpicture}
		\pgfplotstableread{fig-SIR-ER-AvgSusInfRec-All-C.txt}\mainTable
		\begin{axis}[
			ymin=0,ymax=1,xmin=0,xmax=12,
			xlabel={$t$},
			ylabel={Averaged $\Sus{\nodek}$, $\Inf{\nodek}$, $\Rec{\nodek}$},
%			xtick={0,0.05,0.1,0.15,0.2,0.25,0.3},
%			xticklabels={0,0.05,0.1,0.15,0.2,0.25,0.3},
			%	ytick={0},
			%	extra x ticks={0.5,1},
			%	extra y ticks={0.5,1},
			%	extra x tick labels={$\tfrac{1}{2}$, $1$}, 
			%	extra y tick labels={$\tfrac{1}{2}$, $1$}, 
			width={\textwidth},
			height={0.8\textwidth}
			]
			\addplot[color=black,ultra thick] table[x index = 0, y index = 1] from \mainTable;
			\addplot[color=black,ultra thick,dashed] table[x index = 0, y index = 4] from \mainTable;
			\addplot[color=blue,thick] table[x index = 0, y index = 2] from \mainTable;
			\addplot[color=blue,dashed,thick] table[x index = 0, y index = 5] from \mainTable;
			\addplot[color=red,thin] table[x index = 0, y index = 3] from \mainTable;
			\addplot[color=red,dashed,thin] table[x index = 0, y index = 6] from \mainTable;
			%	\foreach \figLineMaker in {2,3,...,6} {
			%		\addplot[color=black] table[x index = 0, y index = \figLineMaker] from \mainTable;
			%	}
			%	\foreach \figLineMaker in {7,8,...,12} {
			%		\addplot[color=black, dashed] table[x index = 0, y index = \figLineMaker] from \mainTable;
			%	}
		\end{axis}
	\end{tikzpicture}
%	\caption{Averaged probabilities over all nodes of nodes being susceptible (black/thick lines), infectious (blue/medium lines), or recovered (red/thin lines) in an Erd\H{o}s--R\'{e}nyi graph according to the rooted-tree approximation (continuous lines) and the average of $10^5$ Gillespie simulations (dashed lines).}
\caption{}
 \label{fig:SIR-ER-AvgSusInf}
\end{subfigure}

\vspace{12pt}
%	\caption{Need to write caption} \label{fig:SIR-ER}
%\end{figure}
%
%
%\begin{figure}
%	\centering

	\begin{subfigure}[b]{0.49\textwidth}
	\centering
	\tikzsetnextfilename{SIR-AlmostTree-SusSelectedNodes}
	\begin{tikzpicture}
		\pgfplotstableread{fig-SIR-AlmostTree-SusSelectedNodes-All-B.txt}\mainTable
		\begin{axis}[
			ymin=0,ymax=1,xmin=0,xmax=6,
			xlabel={$t$},
			ylabel={$\Sus{\nodek}$},
%			xtick={0,0.05,0.1,0.15,0.2,0.25,0.3},
%			xticklabels={0,0.05,0.1,0.15,0.2,0.25,0.3},
			%	ytick={0},
			%	extra x ticks={0.5,1},
			%	extra y ticks={0.5,1},
			%	extra x tick labels={$\tfrac{1}{2}$, $1$}, 
			%	extra y tick labels={$\tfrac{1}{2}$, $1$}, 
			width={\textwidth},
			height={0.8\textwidth}
			]
			\addplot[color=black,mark=o] table[x index = 0, y index = 1] from \mainTable;
			\addplot[color=black,mark=o,dashed,mark options=solid] table[x index = 0, y index = 7] from \mainTable;
			\addplot[color=blue,mark=+] table[x index = 0, y index = 2] from \mainTable;
			\addplot[color=blue,mark=+,dashed,mark options=solid] table[x index = 0, y index = 8] from \mainTable;
			\addplot[color=red,mark=square] table[x index = 0, y index = 4] from \mainTable;
			\addplot[color=red,mark=square,dashed,mark options=solid] table[x index = 0, y index = 10] from \mainTable;
			\addplot[color=olive,mark=triangle] table[x index = 0, y index = 3] from \mainTable;
			\addplot[color=olive,mark=triangle,dashed,mark options=solid] table[x index = 0, y index = 9] from \mainTable;
			%	\foreach \figLineMaker in {2,3,...,6} {
			%		\addplot[color=black] table[x index = 0, y index = \figLineMaker] from \mainTable;
			%	}
			%	\foreach \figLineMaker in {7,8,...,12} {
			%		\addplot[color=black, dashed] table[x index = 0, y index = \figLineMaker] from \mainTable;
			%	}
		\end{axis}
	\end{tikzpicture}
%	\caption{Probabilities of selected nodes (indicated using different colours and mark styles) in an Erd\H{o}s--R\'{e}nyi graph being susceptible at time $t$, according to the rooted-tree approximation (continuous lines) and the average of $10^5$ Gillespie simulations (dashed lines).} 
	\caption{}
	\label{fig:SIR-AlmostTree-SusSelectedNodes}
\end{subfigure}
\begin{subfigure}[b]{0.49\textwidth}
	\centering
	\tikzsetnextfilename{SIR-AlmostTree-AvgSusInf}
	\begin{tikzpicture}
		\pgfplotstableread{fig-SIR-AlmostTree-AvgSusInfRec-All-B.txt}\mainTable
		\begin{axis}[
			ymin=0,ymax=1,xmin=0,xmax=12,
			xlabel={$t$},
			ylabel={Averaged $\Sus{\nodek}$, $\Inf{\nodek}$, $\Rec{\nodek}$},
			%			xtick={0,0.05,0.1,0.15,0.2,0.25,0.3},
			%			xticklabels={0,0.05,0.1,0.15,0.2,0.25,0.3},
			%	ytick={0},
			%	extra x ticks={0.5,1},
			%	extra y ticks={0.5,1},
			%	extra x tick labels={$\tfrac{1}{2}$, $1$}, 
			%	extra y tick labels={$\tfrac{1}{2}$, $1$}, 
			width={\textwidth},
			height={0.8\textwidth}
			]
			\addplot[color=black,ultra thick] table[x index = 0, y index = 1] from \mainTable;
			\addplot[color=black,ultra thick,dashed] table[x index = 0, y index = 4] from \mainTable;
			\addplot[color=blue,thick] table[x index = 0, y index = 2] from \mainTable;
			\addplot[color=blue,dashed,thick] table[x index = 0, y index = 5] from \mainTable;
			\addplot[color=red,thin] table[x index = 0, y index = 3] from \mainTable;
			\addplot[color=red,dashed,thin] table[x index = 0, y index = 6] from \mainTable;
			%	\foreach \figLineMaker in {2,3,...,6} {
			%		\addplot[color=black] table[x index = 0, y index = \figLineMaker] from \mainTable;
			%	}
			%	\foreach \figLineMaker in {7,8,...,12} {
			%		\addplot[color=black, dashed] table[x index = 0, y index = \figLineMaker] from \mainTable;
			%	}
		\end{axis}
	\end{tikzpicture}
%	\caption{Averaged probabilities over all nodes of nodes being susceptible (black/thick lines), infectious (blue/medium lines), or recovered (red/thin lines) in a tree graph with added edges according to the rooted-tree approximation (continuous lines) and the average of $10^5$ Gillespie simulations (dashed lines).} 
\caption{}
\label{fig:SIR-AlmostTree-AvgSusInf}
\end{subfigure}
\caption{\small 
		Comparisions of the rooted-tree approximation in \eqref{E:RootedTreeApprox-SIR} with simulation results from the average of $10^5$ Gillespie algorithm simulations of the full stochastic SIR model. 
		Two different networks are illustrated: subfigures (a) and (b) show results from an Erd{\H o}s--R\'enyi random graph of 100 nodes with probability of connection 0.05; subfigures (c) and (d) show results from a 30-node random tree (generated from a random Pr\"{u}fer sequence) with 10 additional edges added at random.
		Subfigures (a) and (c) show $\Sus{\nodek}$ for four different nodes: results from the rooted-tree approximation are shown as continuous lines and results from Gillespie simulations are shown as dashed lines;
		different nodes are distinguished using different colours and marker styles.
		Subfigures (b) and (d) show $\Sus{\nodek}$ (thick black lines), $\Inf{\nodek}$ (medium thickness blue lines) and $\Rec{\nodek}$ (thin red lines) averaged over all nodes in the network:
		results from the rooted-tree approximation are shown as continuous lines and results from Gillespie simulations are shown as dashed lines.
		Parameters used are $\lamtotn{} = 1$ and $\gamman{} = 0.1$.
		There is a single node that is infectious at $t=0$ and all other nodes are susceptible.
	} \label{fig:SIR-Comparisons} %\label{fig:SIR-AlmostTree}
\end{figure}

\cref{fig:SIR-Comparisons} shows comparisons of the solution of \eqref{E:RootedTreeApprox-SIR} with results obtained from averaging $10^5$ simulations using the Gillespie algorithm.
As previously, \textsc{Matlab} code is available at \url{https://github.com/cameronlhall/rootedtreeapprox}.
Two different networks are shown: an Erd{\H o}s--R\'enyi (ER) random graph (\cref{fig:SIR-ER-SusSelectedNodes,fig:SIR-ER-AvgSusInf}) and a network that is `almost' a tree (\cref{fig:SIR-AlmostTree-SusSelectedNodes,fig:SIR-AlmostTree-AvgSusInf}) in the sense that it was constructed from a random tree by adding some additional edges at random.

If we think of the Gillespie algorithm results as being the `true' solution, we see from \cref{fig:SIR-ER-SusSelectedNodes,fig:SIR-AlmostTree-SusSelectedNodes} that the rooted-tree approximation does indeed give an upper bound on $\Sus{\nodek}$ for each individual node $\nodek$. 
Throughout \cref{fig:SIR-ER-SusSelectedNodes,fig:SIR-AlmostTree-SusSelectedNodes} we see that the rooted-tree approximation deviates from the true solutions by different amounts at different times for different nodes, but the difference is typically substantial as time goes on.
This deviation is observed for the `almost tree' in \cref{fig:SIR-AlmostTree-SusSelectedNodes} as well as for the ER graph in \cref{fig:SIR-ER-SusSelectedNodes}, although we note that the difference between the approximation and the true solution grows faster and becomes larger in the case of the ER graph.

The overall differences between the rooted-tree approximation and the true solution are best seen in \cref{fig:SIR-ER-AvgSusInf,fig:SIR-AlmostTree-AvgSusInf}.
These show $\Sus{\nodek}(t)$, $\Inf{\nodek}(t)$ and $\Rec{\nodek}(t)$ averaged over all nodes in the network.
As may be anticipated from \cref{fig:SIR-ER-SusSelectedNodes,fig:SIR-AlmostTree-SusSelectedNodes}, the rooted-tree approximation gives a overestimate of $\Sus{\nodek}$ (including the equilibrium $\Sus{\nodek}$ as $t \to \infty$) and underestimates the peak in $\Inf{\nodek}$.
Overall, we see that \eqref{E:RootedTreeApprox-SIR} does indeed give bounds on $\Sus{\nodek}$ but that these bounds are not generally very tight.

\section{Rooted-tree approximation for a generalised SEIR model}
\label{S:SEIR}

\subsection{Preliminaries}
\label{S:Prelims-SEIR}

The Susceptible-Exposed-Infectious-Recovered (SEIR) model is a well-established compartment model in the epidemiological literature \cite{BrauerEpidemiology}. 
The SEIR model differs from the SIR model by the introduction of an `exposed' or `latent' state representing individuals that have encountered the disease but are not yet infectious.
Some SEIR models involve multiple classes of exposed state; such models have been analysed mathematically \cite{Bame2008,Diekmann2010,Guo2012} and applied to modelling certain diseases \cite{Cunniffe2012}.

As with the SIR model, the SEIR model has also been extended to networks  \cite{Kang2020,Liu2017,NewmanNetworks,PastorSatorras2015}.
For the most part, network SEIR models in the published literature involve a single exposed state; however, they can be extended to multiple classes of exposed state in an analogous way to compartment models.
Our analysis of SIR models in \cref{S:SIR} can be extended to SEIR models, including in a general setting with arbitrarily many distinct exposed states.
In this section, we replicate our analysis from the previous section but for generalised SEIR models: we construct a node-based approximation of SEIR contagion dynamics that is exact on rooted trees and that yields an upper bound on $\Sus{\nodek}(t)$ on more general graphs.

%REF TO DISCUSSION OF COMPETING I STATES AND OTHER LIMITATIONS
%\todo{Link to later section?}

In our generalised network SEIR model, each node represents an individual, so that at any time a node can either be susceptible (S), exposed of class $\classu$ (E$^{(\classu)}$), infectious (I), or recovered (R).
We assume that there are finitely many ($\numClasses$) different classes of exposed nodes.
Susceptible nodes in contact with infectious nodes may become exposed (in any class) or infectious; 
we refer to the process of a susceptible node changing its state as `infection' regardless of whether the node becomes exposed or infectious.
Exposed nodes may change to a different class of exposed, become infectious, or recover; we assume that exposed nodes cannot become susceptible.
Infectious nodes may recover, but cannot become exposed or susceptible.
Once a node has recovered, it remains recovered for all time.

%The state of a node will change over time based on its current state and the states of its neighbours.
%We assume that all changes of state are Markovian and occur in continuous time ($t$).

%THIS BIT NEEDS TO BE MADE CLEARER AND BETTER TO GET TO THE EQUATIONS. PROBABLY MORE REFERENCES TO THE LITERATURE AS WELL.

Each of these transitions is governed by a different rate parameter. 
The rate of infection (\ie, the total rate at which a susceptible node in contact with an infectious node becomes exposed or infectious) is given by $\lamtotn{}$. 
The probability that a susceptible node becomes exposed of class $\classu$ when infection occurs is given by $\phijn{\classu}{}$; hence, the probability that a susceptible node becomes infectious when infection occurs is $1 - \sum_{\classu} \phijn{\classu}{}$.
The rate at which an exposed node of class $\classu$ becomes an exposed node of class $\classv$ is given by $\ajkn{\classv}{\classu}{}$.
The rate at which an exposed node of class $\classu$ becomes infectious is given by $\mujn{\classu}{}$.
The rate at which an exposed node of class $\classu$ recovers is given by $\nujn{\classu}{}$.
The rate at which an infectious node recovers is given by $\gamman{}$.
These different transitions are summarised below:
%\begin{align*}
%	\text{S (in contact with I)}  &\xrightarrow[\qquad \qquad \quad]{\lamtotn{} \phijn{\classu}{}} \text{E}^{(\classu)} & 
%	\text{S (in contact with I)}  &\xrightarrow[\qquad \qquad \quad]{\lamtotn{} (1 - \sum \phijn{\classu}{})} \text{I} \\
%	\text{E}^{(\classu)} &\xrightarrow[\qquad \qquad \quad]{\ajkn{\classv}{\classu}{}} \text{E}^{(\classv)} &  
%	\text{E}^{(\classu)} &\xrightarrow[\qquad \qquad \quad]{\mujn{\classu}{}} \text{I} \\
%	\text{E}^{(\classu)} &\xrightarrow[\qquad \qquad \quad]{\nujn{\classu}{}} \text{R} &
%	\text{I} &\xrightarrow[\qquad \qquad \quad]{\gamman{}} \text{R} 
%\end{align*}
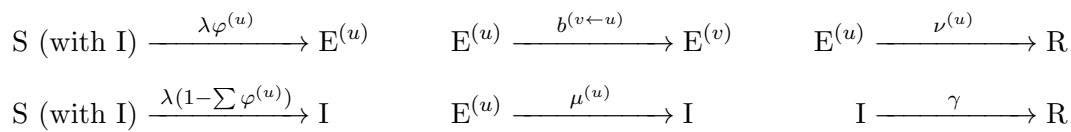
\begin{align*}
	\text{S (with I)}  &\xrightarrow[\qquad \qquad \quad]{\lamtotn{} \phijn{\classu}{}} \text{E}^{(\classu)} & 
	\text{E}^{(\classu)} &\xrightarrow[\qquad \qquad \quad]{\ajkn{\classv}{\classu}{}} \text{E}^{(\classv)} &
	\text{E}^{(\classu)} &\xrightarrow[\qquad \qquad \quad]{\nujn{\classu}{}} \text{R} \\
	\text{S (with I)}  &\xrightarrow[\qquad \qquad \quad]{\lamtotn{} (1 - \sum \phijn{\classu}{})} \text{I} &
	\text{E}^{(\classu)} &\xrightarrow[\qquad \qquad \quad]{\mujn{\classu}{}} \text{I} &
	\text{I} &\xrightarrow[\qquad \qquad \quad]{\gamman{}} \text{R} 
\end{align*}

As in \cref{S:Prelims-SIR}, we assume that the model parameters can depend on the relevant edge or node, and we represent this using subscripts.
The most general approach would be to permit both $\lamtotn{}$ and $\phijn{\classu}{}$ to be edge-dependent; 
however, this level of generality in $\phijn{\classu}{}$ would lead to a problem with the bounding argument in \cref{S:Bounds-SEIR}.
To circumvent this, we permit $\phijn{\classu}{}$ to depend on the recipient node but not on the infecting node; 
that is, we assume $\phijnp{\classu}{\nodek}{\nodej} = \phijn{\classu}{\nodek}$. 
Physically, this would correspond to a situation where 
%different levels of contact and contagiousness mean that the rate of infection, $\lamtotn{}$, depends on the pair of individuals (the directed edge $\nodek \leftarrow \nodej$), but where how 
individual responses to infection (\eg, whether an individual immediately becomes infectious or whether they first enter an exposed state) may vary between individuals but do not depend on the source of infection.

To assist with the analysis of the $\numClasses$ different classes of exposed state, we introduce the $\numClasses$-dimensional vectors $\ExpVec{\nodek}(t)$, $\phivecn{\nodek}$, $\nuvecn{\nodek}$, $\muvecn{\nodek}$, $\uvec$, and $\zerovec$ so that
\begin{align}
	\ExpVec{\nodek}(t) &= 
	\begin{bmatrix}
		\Exp{\nodek}{1}(t) \\
		\Exp{\nodek}{2}(t) \\
		\vdots \\
		\Exp{\nodek}{\numClasses}(t)
	\end{bmatrix},
	&
	\phivecn{\nodek} &=
	\begin{bmatrix}
		\phijn{1}{\nodek} \\
		\phijn{2}{\nodek} \\
		\vdots \\
		\phijn{\numClasses}{\nodek}
	\end{bmatrix},
	&
	\nuvecn{\nodek}&=
	\begin{bmatrix}
		\nujn{1}{\nodek} \\
		\nujn{2}{\nodek} \\
		\vdots \\
		\nujn{\numClasses}{\nodek}
	\end{bmatrix},
	\\
	\muvecn{\nodek}&=
	\begin{bmatrix}
		\mujn{1}{\nodek} \\
		\mujn{2}{\nodek} \\
		\vdots \\
		\mujn{\numClasses}{\nodek}
	\end{bmatrix}, &
	\uvec&=
	\begin{bmatrix}
		1 \\
		1 \\
		\vdots \\
		1
	\end{bmatrix}, &
	\zerovec &=
	\begin{bmatrix}
		0 \\
		0 \\
		\vdots \\
		0
	\end{bmatrix}.
\end{align}
We note that $0 \leq \uvec \cdot \phivecn{\nodek} \leq 1$ for all $\nodek$, and that the rate at which a susceptible node $\nodek$ in contact with an infectious node $\nodej$ becomes infectious is given by
\begin{equation}
	\lamtotnp{\nodek}{\nodej} \left(1 - \sum_{\classu=1}^{\numClasses} \phijn{\classu}{\nodek}\right) = 	\lamtotnp{\nodek}{\nodej} \left(1 - \uvec \cdot \phivecn{\nodek} \right). 
\end{equation}

Lastly, we define the $\numClasses$-by-$\numClasses$ matrix $\An{\nodek}$ so that
\begin{equation}
	\left[\An{\nodek}\right]_{\classu\classv} =
	\begin{cases}
		\mujn{\classv}{\nodek} + \nujn{\classv}{\nodek} + \displaystyle \sum_{\substack{w = 1 \\ w \neq \classv}}^{\numClasses} \ajkn{w}{\classv}{\nodek}, & \classu = \classv, \\
		-\ajkn{\classu}{\classv}{\nodek}; & \classu \neq \classv. \label{E:AnDefn}
	\end{cases}
\end{equation}
%We note that $\An{\nodek}$ is column diagonally dominant with positive values on the diagonal. 
%As a result, the Gershgorin circle theorem guarantees that all eigenvalues of $\An{\nodek}$ have positive real part.

With this notation, the dynamics of contagion on any network can be described using the following equations:
\begin{subequations}
	\label{E:GeneralEvol-SEIR}
\begin{align}
	\diff{\Sus{\nodek}}{t}    &= - \sum_{j \in \Neigh(k)} \lamtotnp{\nodek}{\nodej} \InfSus{\nodej}{\nodek}, \label{E:SusEvol-SEIR}\\
	\diff{\ExpVec{\nodek}}{t} &= \phivecn{\nodek} \sum_{j \in \Neigh(k)} \lamtotnp{\nodek}{\nodej} \InfSus{\nodej}{\nodek} -  \An{\nodek} \ExpVec{\nodek}, \label{E:ExpEvol-SEIR} \\
	\diff{\Inf{\nodek}}{t}    &= (1 - \uvec \cdot \phivecn{\nodek}) \sum_{j \in \Neigh(k)} \lamtotnp{\nodek}{\nodej} \InfSus{\nodej}{\nodek} + \muvecn{\nodek} \cdot \ExpVec{\nodek} - \gamman{\nodek} \Inf{\nodek},  \label{E:InfEvol-SEIR}\\
	\diff{\Rec{\nodek}}{t}    &= \nuvecn{\nodek} \cdot \ExpVec{\nodek} + \gamman{\nodek} \Inf{\nodek}, \label{E:RecEvol-SEIR}
\end{align}
\end{subequations}
which must be solved subject to suitable initial conditions.

Note that if $\phivecn{}$ were permitted to depend on the source of infection as well as on the node that becomes infected then the corresponding $\phivecn{\nodek \leftarrow \nodej}$ terms would need to be included inside the summations in equations \eqref{E:ExpEvol-SEIR} and \eqref{E:InfEvol-SEIR}. 

Note also that \eqref{E:SusEvol-SEIR} can be used to express \eqref{E:ExpEvol-SEIR} and \eqref{E:InfEvol-SEIR} in the equivalent forms
\begin{subequations}
	\setcounter{equation}{1}
	\label{E:GeneralEvol-SEIRAlt}
	\begin{align}
		\diff{\ExpVec{\nodek}}{t} &= -\phivecn{\nodek} \diff{\Sus{\nodek}}{t} -  \An{\nodek} \ExpVec{\nodek}, \label{E:ExpEvol-SEIRAlt} \\
		\diff{\Inf{\nodek}}{t}    &= -(1 - \uvec \cdot \phivecn{\nodek}) \diff{\Sus{\nodek}}{t} + \muvecn{\nodek} \cdot \ExpVec{\nodek} - \gamman{\nodek} \Inf{\nodek}.  \label{E:InfEvol-SEIRAlt}
	\end{align}
\end{subequations}
Given the length of the expressions that we obtain for $\lamtotnp{\nodek}{\nodej} \InfSus{\nodej}{\nodek}$ in our analysis, we will sometimes prefer \eqref{E:ExpEvol-SEIRAlt} and \eqref{E:InfEvol-SEIRAlt} over \eqref{E:ExpEvol-SEIR} and \eqref{E:InfEvol-SEIR} for concision.

\subsection{Exact SEIR dynamics on a rooted tree}
\label{S:RootedTree-SEIR}

As in \cref{S:RootedTree-SIR}, we begin by considering contagion dynamics on a rooted tree, where there is a single node, $k = 0$, which is the source of infection. 
This node may either be exposed or infectious at $t = 0$.
Introducing equivalent notation and following the same logic as for the derivation of \eqref{E:RootedEvol-SIR}, we find that the evolution equations for node state probabilities on a rooted tree are
\begin{subequations}
	\label{E:RootedEvol-SEIR}
\begin{align}
	\diff{\Sus{\nodek}}{t}    &=
	\begin{cases}
		0, & \nodek = 0, \\
		-\lamtotn{\nodek} \InfSus{\parn{\nodek}}{\nodek}, & \nodek \neq 0; 
	\end{cases} \label{E:SusEvol-SEIR-Rooted} \\
	\diff{\ExpVec{\nodek}}{t} &= 
	\begin{cases}
		- \An{\nodek} \ExpVec{\nodek},  & \nodek = 0, \\
		\lamtotn{\nodek} \phivecn{\nodek} \InfSus{\parn{\nodek}}{\nodek} - \An{\nodek} \ExpVec{\nodek},  & \nodek \neq 0;
	\end{cases} \label{E:ExpEvol-SEIR-Rooted} \\
	\diff{\Inf{\nodek}}{t}    &=
	\begin{cases}
		\muvecn{\nodek} \cdot \ExpVec{\nodek} - \gamman{\nodek} \Inf{\nodek},  & \nodek = 0, \\
		\lamtotn{\nodek}(1 - \uvec \cdot \phivecn{\nodek}) \InfSus{\parn{\nodek}}{\nodek} + \muvecn{\nodek} \cdot \ExpVec{\nodek} - \gamman{\nodek} \Inf{\nodek},  & \nodek \neq 0 ;
	\end{cases}  \label{E:InfEvol-SEIR-Rooted}\\
	\diff{\Rec{\nodek}}{t}    &= \nuvecn{\nodek} \cdot \ExpVec{\nodek} + \gamman{\nodek} \Inf{\nodek}. \label{E:RecEvol-SEIR-Rooted}
\end{align}
\end{subequations}

These equations need to be solved subject to initial conditions where
\begin{align}
	\Sus{\nodek}(0) &= 1, & \ExpVec{\nodek}(0) &=  \zerovec, & \Inf{\nodek}(0) &= \Rec{\nodek}(0) = 0, & \nodek &\neq 0, \label{E:MainICs}
\end{align}
and where $\ExpVec{0}(0) = \ExpVecInit$ and $\Inf{0}(0) = \InfInit$ are specified, but $\Sus{0}(0) = \Rec{0}(0) = 0$.
%Since it is impossible for infection to spread to the root from any other node, the evolution of the node state probabilities at $\nodek = 0$ are given by
%\begin{align}
%	\diff{\Sus{0}}{t}    &= 0, \label{E:Sus0Evol}\\
%	\diff{\ExpVec{0}}{t} &=  - \An{0} \ExpVec{0}, \label{E:Exp0Evol}\\
%	\diff{\Inf{0}}{t}    &= \muvecn{0} \cdot \ExpVec{0} - \gamman{0} \Inf{0}, \label{E:Inf0Evol}\\
%	\diff{\Rec{0}}{t}    &= \nuvecn{0} \cdot \ExpVec{0} + \gamman{0} \Inf{0}. \label{E:Rec0Evol}
%\end{align}
%For the root, we will only consider initial conditions where $\Sus{0}(0) = 0$ and $\Rec{0}(0) = 0$; the results for any other scenario can be obtained by a simple rescaling. Thus, the general initial conditions for the root node are
%\begin{align}
%	\Sus{0}(0) &= \Rec{0}(0) = 0, & \ExpVec{0}(0) &=  \ExpVecInit, & \Inf{0}(0) &= \InfInit. \label{E:n=0ICs}
%\end{align}
We note that $\phivecn{0}$ does not appear in system \eqref{E:RootedEvol-SEIR} or in the initial conditions. 
As we will see, it will be convenient to define $\phivecn{0}$ so that $\phivecn{0} = \ExpVecInit$, and hence $\InfInit = 1 - \uvec \cdot \phivecn{0}$.  

System \eqref{E:RootedEvol-SEIR} is not closed because of the presence of $\InfSus{\parn{\nodek}}{\nodek}$.
As in \cref{S:RootedTree-SIR}, we exploit the properties of a rooted tree to find an expression for $\InfSus{\parn{\nodek}}{\nodek}$ in terms of the node state probabilities and hence obtain a closed system. Since the parent node of node 0 is not defined, we assume (unless otherwise specified) that $\nodek \neq 0$ in all analysis below where $\parn{\nodek}$ is mentioned.

We begin by noting that the law of total probability gives
\begin{equation}
\Sus{\nodek} = \SusSus{\parn{\nodek}}{\nodek} + \sum_{\classu=1}^{\numClasses} \ExpSus{\parn{\nodek}}{\classu}{\nodek} + \InfSus{\parn{\nodek}}{\nodek} + \RecSus{\parn{\nodek}}{\nodek}. \label{E:SusExhaustionSEIR}
\end{equation}
The fact that infection can only spread from node $\parn{\nodek}$ to node $\nodek$ and not \textit{vice versa} means that if either $X_{\parn{\nodek}} = \text{S}$ or $X_{\parn{\nodek}} = \text{E}^{(j)}$ then $X_{\nodek} = \text{S}$.
Hence, $\SusSus{\parn{\nodek}}{\nodek} = \Sus{\parn{\nodek}}$ and $\ExpSus{\parn{\nodek}}{\classu}{\nodek} = \Exp{\parn{\nodek}}{\classu}$.
Thus, \eqref{E:SusExhaustionSEIR} can be rearranged to give
\begin{align}
\InfSus{\parn{\nodek}}{\nodek}(t) &= \Sus{\nodek} - \Sus{\parn{\nodek}}(t) - \sum_{\classu=1}^{\numClasses} \Exp{\parn{\nodek}}{\classu} -  \RecSus{\parn{\nodek}}{\nodek} \notag \\
&= \Sus{\nodek} - \Sus{\parn{\nodek}} - \uvec \cdot \ExpVec{\parn{\nodek}} -  \RecSus{\parn{\nodek}}{\nodek}.  \label{E:InfSusFromRecSus}
\end{align}

As previously, we now seek a differential equation for $\RecSus{\parn{\nodek}}{\nodek}$ that can be directly integrated to obtain  $\RecSus{\parn{\nodek}}{\nodek}$ in terms of node state probabilities.
The only way to achieve a state where $X_{\parn{\nodek}} = \text{R}$ and $X_{\nodek} = \text{S}$ is for node $\parn{\nodek}$ to recover (either form an exposed state or an infectious state) while node $\nodek$ is susceptible. 
Once node $\parn{\nodek}$ has recovered, this state will then be permanent.
Since $\ExpSus{\parn{\nodek}}{\classu}{\nodek} = \Exp{\parn{\nodek}}{\classu}$, it therefore follows that
\begin{align}
\diff{\RecSus{\parn{\nodek}}{\nodek}}{t}    &= \nuvecn{\parn{\nodek}} \cdot \ExpVec{\parn{\nodek}} + \gamman{\parn{\nodek}} \InfSus{\parn{\nodek}}{\nodek}. 
\end{align}
Using \eqref{E:SusEvol-SEIR-Rooted}, this rearranges to give
\begin{align}
\diff{\RecSus{\parn{\nodek}}{\nodek}}{t}    &= \nuvecn{\parn{\nodek}} \cdot \ExpVec{\parn{\nodek}} - \frac{\gamman{\parn{\nodek}}}{\lamtotn{\nodek}} \diff{\Sus{\nodek}}{t}, \label{E:RecSusEvol-SEIR}
\end{align}

The next step is to rewrite $\nuvecn{\parn{\nodek}} \cdot \ExpVec{\parn{\nodek}}(t)$ in terms of the derivatives of node state probabilities.
For any node $\nodek$ (including $\nodek = 0$), let $\Mn{\nodek}$ be the block matrix defined by
\begin{equation}
\Mn{\nodek} =
\begin{bmatrix}
1 & \zerovec^T \\
- \phivecn{\nodek} & \An{\nodek}
\end{bmatrix},
\label{E:MDefn}
\end{equation}
so that the block matrix inversion formula \cite{MatrixCookbook} gives
\begin{equation}
\Mn{\nodek}^{-1} =
\begin{bmatrix}
1 & \zerovec^T \\[12pt]
\An{\nodek}^{-1} \phivecn{\nodek} & \An{\nodek}^{-1}
\end{bmatrix}.
\end{equation}
Using $\Mn{\nodek}$, we can rewrite equations \eqref{E:SusEvol-SEIR-Rooted} and \eqref{E:ExpEvol-SEIR-Rooted} together as
\begin{equation}
\begin{bmatrix}
\diff{\Sus{\nodek}}{t} \\[6pt]
\diff{\ExpVec{\nodek}}{t}
\end{bmatrix}
=
- \Mn{\nodek}
\begin{bmatrix}
\lamtotn{\nodek} \InfSus{\parn{\nodek}}{\nodek} \\[6pt]
\ExpVec{\nodek}
\end{bmatrix}. \label{E:SusAndExpEvol}
\end{equation}
If we assert that $\InfSus{p(0)}{0}(t) \equiv 0$, then %comparison with the $\nodek = 0$ case in equations \eqref{E:SusEvol-SEIR-Rooted} and \eqref{E:ExpEvol-SEIR-Rooted} shows that 
\eqref{E:SusAndExpEvol} also applies when $\nodek = 0$.

We now use $\Mn{\parn{\nodek}}$ to express  $\nuvecn{\parn{\nodek}} \cdot \ExpVec{\parn{\nodek}}$ in terms of derivatives as follows:
\begin{align}
\nuvecn{\parn{\nodek}} \cdot \ExpVec{\parn{\nodek}}
&= 
\begin{bmatrix}
0 & \nuvecn{\parn{\nodek}}^T
\end{bmatrix}
\begin{bmatrix}
\lamtotn{\parn{\nodek}}\InfSus{\pparn{\nodek}}{\parn{\nodek}} \\[6pt]
\ExpVec{\parn{\nodek}}
\end{bmatrix} \label{E:nuDotExpTemp}\\[6pt]
&= 
- \begin{bmatrix}
0 & \nuvecn{\parn{\nodek}}^T
\end{bmatrix}
\Mn{\parn{\nodek}}^{-1}
\begin{bmatrix}
\diff{\Sus{\parn{\nodek}}}{t} \\[6pt]
\diff{\ExpVec{\parn{\nodek}}}{t}
\end{bmatrix} \notag \\[6pt]
%&= - 	\begin{bmatrix}
%	0 & \nuvecn{\parn{\nodek}}^T
%\end{bmatrix}
%	\begin{bmatrix}
%	\frac{1}{\lamtotn{\parn{\nodek}}} & \zerovec^T \\[12pt]
%	\frac{\An{\parn{\nodek}}^{-1} \lamvecn{\parn{\nodek}}}{\lamtotn{\parn{\nodek}}} & \An{\parn{\nodek}}^{-1}
%\end{bmatrix}
%	\begin{bmatrix}
%	\diff{\Sus{\parn{\nodek}}}{t} \\[6pt]
%	\diff{\ExpVec{\parn{\nodek}}}{t}
%\end{bmatrix} \\[6pt]
%&= - 	\begin{bmatrix}
%	\frac{ \nuvecn{\parn{\nodek}}^T \An{\parn{\nodek}}^{-1} \lamvecn{\parn{\nodek}}}{\lamtotn{\parn{\nodek}}} & \nuvecn{\parn{\nodek}}^T \An{\parn{\nodek}}^{-1}
%\end{bmatrix}
%	\begin{bmatrix}
%	\diff{\Sus{\parn{\nodek}}}{t} \\[6pt]
%	\diff{\ExpVec{\parn{\nodek}}}{t}
%\end{bmatrix} \\[6pt]
&= -\nuvecn{\parn{\nodek}}^T \An{\parn{\nodek}}^{-1} \phivecn{\parn{\nodek}} \diff{\Sus{\parn{\nodek}}}{t} - \nuvecn{\parn{\nodek}}^T \An{\parn{\nodek}}^{-1} \diff{\ExpVec{\parn{\nodek}}}{t}. \label{E:nuDotExpTemp2}
\end{align}
Note that equation \eqref{E:nuDotExpTemp2} applies even when $\parn{\nodek} = 0$; even though the value of $\InfSus{\pparn{\nodek}}{\parn{\nodek}}$ would be undefined in \eqref{E:nuDotExpTemp}, it is multiplied by zero and does not affect the final result.
%Additionally, equation \eqref{E:nuDotExpTemp2} is valid when $\parn{\nodek} = 0$ since $\diff{\Sus{0}}{t} = 0$.

Substituting \eqref{E:nuDotExpTemp2} into \eqref{E:RecSusEvol-SEIR} yields
\begin{equation}
\diff{\RecSus{\parn{\nodek}}{\nodek}}{t}   = -\nuvecn{\parn{\nodek}}^T \An{\parn{\nodek}}^{-1} \phivecn{\parn{\nodek}} \diff{\Sus{\parn{\nodek}}}{t} - \nuvecn{\parn{\nodek}}^T \An{\parn{\nodek}}^{-1} \diff{\ExpVec{\parn{\nodek}}}{t} - \frac{\gamman{\parn{\nodek}}}{\lamtotn{\nodek}} \diff{\Sus{\nodek}}{t},
\end{equation}
and hence we find that 
\begin{equation}
\RecSus{\parn{\nodek}}{\nodek}   = \Cn{\nodek} -\nuvecn{\parn{\nodek}}^T \An{\parn{\nodek}}^{-1} \phivecn{\parn{\nodek}} \Sus{\parn{\nodek}} - \nuvecn{\parn{\nodek}}^T \An{\parn{\nodek}}^{-1} \ExpVec{\parn{\nodek}} - \frac{\gamman{\parn{\nodek}}}{\lamtotn{\nodek}} \Sus{\nodek},
\label{E:RecSusTemp1}
\end{equation}
where $\Cn{\nodek}$ is a constant to be determined from the initial conditions.

In the case where $\parn{\nodek} \neq 0$, the initial conditions in \eqref{E:MainICs} yield
\begin{equation}
\Cn{\nodek} =  \nuvecn{\parn{\nodek}}^T \An{\parn{\nodek}}^{-1} \phivecn{\parn{\nodek}} + \frac{\gamman{\parn{\nodek}}}{\lamtotn{\nodek}}. \label{E:CnFormula}
\end{equation}
In the case where $\parn{\nodek} = 0$, the initial conditions yield
\begin{equation}
\Cn{\nodek} =  \nuvecn{0}^T \An{0}^{-1} \ExpVecInit + \frac{\gamman{0}}{\lamtotn{\nodek}}.
\end{equation}
As noted previously, this motivates us to define $\phivecn{0} = \ExpVecInit$ so that \eqref{E:CnFormula} can be used to give the constant $\Cn{\nodek}$ for all nodes $\nodek \neq 0$.

%Since node $0$ does not have a parent, the transmission parameters $\lamvecn{0}$, $\lamstarn{0}$, and $\lamtotn{0}$ are not defined.
%If we assert that these transmission parameters are chosen so that
%\begin{equation}
%\frac{\lamvecn{0}}{\lamtotn{0}} = \ExpVecInit, 
%\end{equation}
%then equation \eqref{E:CnFormula} can be used to define the parameter $\Cn{\nodek}$ for all $n \neq 0$.
%As we will see, this affects our choice of initial conditions in the case where we assume that the model parameters do not vary between nodes.

Combining \eqref{E:RecSusTemp1} and \eqref{E:CnFormula}, we obtain
%\begin{multline}
%\RecSus{\parn{\nodek}}{\nodek} =
%\left[\An{\parn{\nodek}}^{-T} \nuvecn{\parn{\nodek}} \right] \cdot \phivecn{\parn{\nodek}} 
%+ \frac{\gamman{\parn{\nodek}}}{\lamtotn{\nodek}} \\
%-\left[\An{\parn{\nodek}}^{-T} \nuvecn{\parn{\nodek}} \right] \cdot \phivecn{\parn{\nodek}} \Sus{\parn{\nodek}}
%- \left[\An{\parn{\nodek}}^{-T} \nuvecn{\parn{\nodek}} \right] \cdot \ExpVec{\parn{\nodek}}
%- \frac{\gamman{\parn{\nodek}}}{\lamtotn{\nodek}} \Sus{\nodek}. \label{E:RecSusTemp2}
%\end{multline}
an expression for $\RecSus{\parn{\nodek}}{\nodek}$ that can be substituted into \eqref{E:InfSusFromRecSus} to yield
%SAD TIMES
%\begin{multline}
%	\InfSus{\parn{\nodek}}{\nodek} = 
%	- \left[\An{\parn{\nodek}}^{-T} \nuvecn{\parn{\nodek}} \right] \cdot \phivecn{\parn{\nodek}} 
%	- \frac{\gamman{\parn{\nodek}}}{\lamtotn{\nodek}} 
%	+ \frac{\lamtotn{\nodek} + \gamman{\parn{\nodek}}}{\lamtotn{\nodek}}\Sus{\nodek} \\
%	- \left(1 - \left[\An{\parn{\nodek}}^{-T} \nuvecn{\parn{\nodek}} \right] \cdot \phivecn{\parn{\nodek}} \right)\Sus{\parn{\nodek}} 
%	- \left( \uvec  - \An{\parn{\nodek}}^{-T} \nuvecn{\parn{\nodek}} \right) \cdot \ExpVec{\parn{\nodek}}. 
%\end{multline}
\begin{multline}
\InfSus{\parn{\nodek}}{\nodek} = 
- \nuvecn{\parn{\nodek}}^T\An{\parn{\nodek}}^{-1} \phivecn{\parn{\nodek}} 
- \frac{\gamman{\parn{\nodek}}}{\lamtotn{\nodek}} 
+ \frac{\lamtotn{\nodek} + \gamman{\parn{\nodek}}}{\lamtotn{\nodek}}\Sus{\nodek} \\
- \left(1 - \nuvecn{\parn{\nodek}}^{T} \An{\parn{\nodek}}^{-1} \phivecn{\parn{\nodek}} \right)\Sus{\parn{\nodek}} 
 - \left( \uvec  - \An{\parn{\nodek}}^{-T} \nuvecn{\parn{\nodek}} \right) \cdot \ExpVec{\parn{\nodek}}. 
 \label{E:InfSusTemp1}
\end{multline}
%As required, this gives us an equation for $\InfSus{\parn{\nodek}}{\nodek}$ that depends only on the node state probabilities.

%For the bounding argument in \cref{S:Bounds-SEIR} it will be important to establish the signs of the coefficients of $\Sus{\parn{\nodek}}$ and $\ExpVec{\parn{\nodek}}$ in \eqref{E:InfSusTemp1}.

We note that \eqref{E:AnDefn} implies that
\begin{equation}
\sum_{\classu=1}^{\numClasses} \left[\An{\nodek}\right]_{\classu\classv} = \mujn{\classv}{\nodek} + \nujn{\classv}{\nodek},
\end{equation}
and hence
%\begin{equation}
$\An{\nodek}^T \uvec = \muvecn{\nodek} + \nuvecn{\nodek}$.  This rearranges to yield 
%\end{equation}
%Rearranging, this gives
$
\An{\nodek}^{-T} \muvecn{\nodek} = \uvec - \An{\nodek}^{-T} \nuvecn{\nodek} %\label{E:AInvTrnu}
$
so that \eqref{E:InfSusTemp1} becomes
\begin{multline}
	\InfSus{\parn{\nodek}}{\nodek} = 
	- \nuvecn{\parn{\nodek}}^T\An{\parn{\nodek}}^{-1} \phivecn{\parn{\nodek}} 
	- \frac{\gamman{\parn{\nodek}}}{\lamtotn{\nodek}} 
	+ \frac{\lamtotn{\nodek} + \gamman{\parn{\nodek}}}{\lamtotn{\nodek}}\Sus{\nodek} \\
	- \left(1 - \nuvecn{\parn{\nodek}}^{T} \An{\parn{\nodek}}^{-1} \phivecn{\parn{\nodek}} \right)\Sus{\parn{\nodek}} 
	- \muvecn{\parn{\nodek}}^T\An{\parn{\nodek}}^{-1} \ExpVec{\parn{\nodek}}. 
	\label{E:InfSusTemp2}
\end{multline}

%and hence
%\begin{align}
%\left[\An{\nodek}^{-T} \nuvecn{\nodek} \right] \cdot \phivecn{\nodek} 
%&= \uvec \cdot \phivecn{\nodek} - \left[\An{\nodek}^{-T} \muvecn{\nodek}\right] \cdot \phivecn{\nodek}. \label{E:AInvTrnuDotlam}
%\end{align}

As an aside, we note from \eqref{E:AnDefn} that $\An{\nodek}^T$ is a strictly diagonally dominant matrix with positive diagonal entries. From \cite{BermanNonnegativeMatrices}, it follows that $\An{\nodek}^T$ is inverse-positive. 
Hence, the elements of $\An{\nodek}^{-T} \muvecn{\nodek}$ and $\An{\nodek}^{-T} \nuvecn{\nodek}$ are all between 0 and 1 (inclusive) and we note that the coefficients of $\Sus{\parn{\nodek}}$ and $\Exp{\parn{\nodek}}{\classu}$ in \eqref{E:InfSusTemp2} are all nonpositive.

Using \eqref{E:InfSusTemp2} and \eqref{E:GeneralEvol-SEIRAlt}, system \eqref{E:RootedEvol-SEIR} can be rearranged to give
\begin{subequations}
	\label{E:RootedEvol-SEIR-Full}
\begin{align}
	%	\diff{\Sus{\nodek}}{t} &+ \left(\lamtotn{\nodek} + \gamman{\parn{\nodek}}\right) \Sus{\nodek}(t) \\
	%	& \qquad = 
	%	\begin{cases}
	%		\\
	%	0, & n = 0, \\
	%	\\
	%	\displaystyle \frac{\lamtotn{\nodek}}{\lamtotn{\parn{\nodek}}} 
	%	\left( \lamstarn{\parn{\nodek}} + \mvecn{\parn{\nodek}} \cdot \lamvecn{\parn{\nodek}} \right)
	%	\Sus{\parn{\nodek}}(t)
	%	%\\ \displaystyle \qquad 
	%	+ \lamtotn{\nodek} \mvecn{\parn{\nodek}} \cdot \ExpVec{\parn{\nodek}}(t)
	%	\\ \displaystyle \qquad 
	%	+ \lamtotn{\nodek} + \gamman{\parn{\nodek}} 
	%	%\\ \displaystyle \qquad 
	%	- \frac{\lamtotn{\nodek}}{\lamtotn{\parn{\nodek}}} 
	%	\left( \lamstarn{\parn{\nodek}} + \mvecn{\parn{\nodek}} \cdot \lamvecn{\parn{\nodek}} \right),
	%	& n \neq 0;
	%	\end{cases}
	\diff{\Sus{\nodek}}{t} &=0, 
	& \nodek&=0, \label{E:Sus0EvolTreeFull} \\
	\diff{\Sus{\nodek}}{t}  &= - \left(\lamtotn{\nodek} + \gamman{\parn{\nodek}}\right) \Sus{\nodek}
	+\lamtotn{\nodek} \left(1 - \nuvecn{\parn{\nodek}}^{T} \An{\parn{\nodek}}^{-1} \phivecn{\parn{\nodek}} \right)\Sus{\parn{\nodek}} \notag \\
	&\qquad + \lamtotn{\nodek} \muvecn{\parn{\nodek}}^T\An{\parn{\nodek}}^{-1} \ExpVec{\parn{\nodek}} 
	+ \gamman{\parn{\nodek}} %\notag \\
	%&\qquad 
	- \lamtotn{\nodek} \nuvecn{\parn{\nodek}}^T\An{\parn{\nodek}}^{-1} \phivecn{\parn{\nodek}}, & \nodek&\neq0, \label{E:SusEvol-SEIR-TreeFull} \\
	\diff{\ExpVec{\nodek}}{t}  &=  -\phivecn{\nodek} \diff{\Sus{\nodek}}{t} - \An{\nodek} \ExpVec{\nodek}, \label{E:ExpEvol-SEIR-TreeFull} \\
	\diff{\Inf{\nodek}}{t} 
	&= 
	-(1 - \uvec\cdot \phivecn{\nodek}) \diff{\Sus{\nodek}}{t} - \gamman{\nodek} \Inf{\nodek}  + \muvecn{\nodek} \cdot \ExpVec{\nodek}, \label{E:InfEvol-SEIR-TreeFull} \\
	\diff{\Rec{\nodek}}{t}
	&= 
	\nuvecn{\nodek} \cdot \ExpVec{\nodek} + \gamman{\nodek} \Inf{\nodek}. \label{E:RecEvol-SEIR-Tree}
\end{align}
\end{subequations}

As for the SIR model in \cref{S:RootedTree-SIR}, this is a partially-decoupled system.
To see this, we observe that the dynamics of $\Sus{\nodek}$ in  \eqref{E:SusEvol-SEIR-TreeFull} are independent of $\ExpVec{\nodek}$; instead, $\diff{\Sus{\nodek}}{t}$ depends only on $\Sus{\nodek}$ and the node state probabilities at the parent node.
Since equations \eqref{E:SusEvol-SEIR-TreeFull} and \eqref{E:ExpEvol-SEIR-TreeFull} are both independent of $\Inf{\nodek}(t)$ and $\Rec{\nodek}(t)$, this implies that \eqref{E:RootedEvol-SEIR-Full} can be solved from the root outwards, with $\Sus{\nodek}$ solved before $\ExpVec{\nodek}$ at each subsequent node.

Moreover, consider the case where exposed states are traversed in order---that is, where $\ajkn{\classu}{\classv}{\nodek}$ is zero whenever $\classu < \classv$).
This situation is physically plausible, since it corresponds to a case where a diseased individual can progress through different exposed ``stages'' before becoming infectious or recovering, but can never return to an earlier class of exposed state from a more advanced class.
In this case, the matrix $\An{\nodek}$ will be lower triangular and hence the scalar equations that constitute \eqref{E:ExpEvol-SEIR-TreeFull} will also be partially decoupled.
Since system \eqref{E:RootedEvol-SEIR-Full} is linear, this implies that the full solution can be obtained exactly by the sequential solving of linear scalar ordinary differential equations;
it is not even necessary to solve an eigenvalue problem in order to obtain the exact solution to SEIR dynamics on a rooted tree.
While we do not present closed-form solutions here, it is theoretically possible to obtain results analogous to \eqref{E:RootedEvol-SEIR-Full} using standard methods for nonhomogeneous constant-coefficients differential equations.

As in \cref{S:Solutions-SIR}, we test the rooted-tree formulation in system \eqref{E:RootedEvol-SEIR-Full} by considering SEIR dynamicson a chain.
For simplicity, we consider the case where there is a single class of exposed state and so the vectors and matrices in \eqref{E:RootedEvol-SEIR-Full} can be replaced by scalars.
Noting that the equivalent of $\An{\nodek}$ will be $\musimpn{\nodek} + \nusimpn{\nodek}$, this leads to the system
\begin{subequations}
	\label{E:ChainEvol-SEIR-Simp}
	\begin{align}
		\diff{\Sus{\nodek}}{t} &=0, 
		& \nodek&=0, \label{E:Sus0EvolTreeSimp} \\
		\diff{\Sus{\nodek}}{t}  &= 
		\frac{\lamtotn{\nodek}  \musimpn{\nodek-1}  }{\musimpn{\nodek-1} + \nusimpn{\nodek-1}} \left(\phin{\nodek-1}\Sus{\nodek-1} + \ExpSimp{\nodek-1} \right) \notag \\ 
		&\qquad - \left(\lamtotn{\nodek} + \gamman{\nodek-1}\right) \Sus{\nodek} + \gamman{\nodek-1} %\notag \\
		%&\qquad 
		- \frac{\lamtotn{\nodek}  \nusimpn{\nodek-1} \phin{\nodek-1} }{\musimpn{\nodek-1} + \nusimpn{\nodek-1}}, & \nodek&\neq0, \label{E:SusEvol-SEIR-TreeSimp} \\
		\diff{\ExpSimp{\nodek}}{t}  &=  -\phin{\nodek} \diff{\Sus{\nodek}}{t} - (\musimpn{\nodek} + \nusimpn{\nodek}) \ExpSimp{\nodek}, \label{E:ExpEvol-SEIR-TreeSimp} \\
		\diff{\Inf{\nodek}}{t} 
		&= 
		-(1 - \phin{\nodek}) \diff{\Sus{\nodek}}{t} - \gamman{\nodek} \Inf{\nodek}  + \musimpn{\nodek} \ExpSimp{\nodek}, \label{E:InfEvol-SEIR-TreeSimp} %\\
%		\diff{\Rec{\nodek}}{t}
%		&= 
%		\nusimpn{\nodek} \ExpSimp{\nodek} + \gamman{\nodek} \Inf{\nodek}. \label{E:RecEvol-SEIR-TreeSimp}
	\end{align}
\end{subequations}

\cref{fig:SEIR-Chain-SusInfAllNodes} shows a comparison of $\Sus{\nodek}(t)$ and $\Inf{\nodek}(t)$ obtained from the numerical solution of \eqref{E:ChainEvol-SEIR-Simp} with the average of $10^5$ Gillespie algorithm simulations of the underlying stochastic model
(code again available at \url{https://github.com/cameronlhall/rootedtreeapprox}).
As in \cref{fig:SIR-Chain-SusInfAllNodes}, this exemplifies the fact that system \eqref{E:ChainEvol-SEIR-Simp} is exact; the two sets of results are virtually indistinguishable.

\begin{figure}
	\centering
	\begin{subfigure}{0.49\textwidth}
		\centering
		\tikzsetnextfilename{SEIR-Chain-SusAllNodes}
		\begin{tikzpicture}
			\pgfplotstableread{fig-SEIR-Chain-SusAllNodes-DE-B.txt}\deSolutionTable
			\pgfplotstableread{fig-SEIR-Chain-SusAllNodes-Gil-B.txt}\gilSolutionTable
			\begin{axis}[
				ymin=0,ymax=1,xmin=0,xmax=15,
				xlabel={$t$},
				ylabel={$\Sus{\nodek}$},
				%	xtick={0},
				%	ytick={0},
				%	extra x ticks={0.5,1},
				%	extra y ticks={0.5,1},
				%	extra x tick labels={$\tfrac{1}{2}$, $1$}, 
				%	extra y tick labels={$\tfrac{1}{2}$, $1$}, 
				width={\textwidth},
				height={0.8\textwidth}
				]
				\addplot[color=blue, very thick] table[x index = 0, y index = 2] from \deSolutionTable;
				\addplot[color=blue, mark=+, only marks] table[x index = 0, y index = 2] from \gilSolutionTable;
				\foreach \figLineMaker in {1,3,4,5,...,11} {
					\addplot[color=black] table[x index = 0, y index = \figLineMaker] from \deSolutionTable;
					\addplot[color=black, mark=+, only marks] table[x index = 0, y index = \figLineMaker] from \gilSolutionTable;
				}
			\end{axis}
		\end{tikzpicture}
		%		\caption{Probability of nodes being susceptible at time $t$.}
		\caption{}
		\label{fig:SEIR-Chain-SusAllNodes}
	\end{subfigure}%
	\hspace{\fill}
	\begin{subfigure}{0.49\textwidth}
		\centering
		\tikzsetnextfilename{SEIR-Chain-InfAllNodes}
		\begin{tikzpicture}
			\pgfplotstableread{fig-SEIR-Chain-InfAllNodes-DE-B.txt}\deSolutionTable
			\pgfplotstableread{fig-SEIR-Chain-InfAllNodes-Gil-B.txt}\gilSolutionTable
			\begin{axis}[
				ymin=0,ymax=0.8,xmin=0,xmax=30,
				xlabel={$t$},
				ylabel={$\Inf{\nodek}$},
				%	xtick={0},
				%	ytick={0},
				%	extra x ticks={0.5,1},
				%	extra y ticks={0.5,1},
				%	extra x tick labels={$\tfrac{1}{2}$, $1$}, 
				%	extra y tick labels={$\tfrac{1}{2}$, $1$}, 
				width={\textwidth},
				height={0.8\textwidth}
				]
				\addplot[color=blue, very thick] table[x index = 0, y index = 2] from \deSolutionTable;
				\addplot[color=blue, mark=+, only marks] table[x index = 0, y index = 2] from \gilSolutionTable;
				\foreach \figLineMaker in {1,3,4,5,...,11} {
					\addplot[color=black] table[x index = 0, y index = \figLineMaker] from \deSolutionTable;
					\addplot[color=black, mark=+, only marks] table[x index = 0, y index = \figLineMaker] from \gilSolutionTable;
				}
			\end{axis}
		\end{tikzpicture}
		%		\caption{Probability of nodes being infectious at time $t$.} 
		\caption{}
		\label{fig:SEIR-Chain-InfAllNodes}
		%\label{TK}
	\end{subfigure}
	\caption{\small 
		Comparision of the rooted-tree solutions for $\Sus{\nodek}$ and $\Inf{\nodek}$ based on numerical solution of  \eqref{E:ChainEvol-SEIR-Simp} with simulation results from the average of $10^5$ Gillespie algorithm simulations of the full stochastic model. 
		Subfigure (a) shows results for $\Sus{\nodek}$ while subfigure (b) shows results for $\Inf{\nodek}$.
		In both cases, the rooted tree solutions are shown as continuous lines and the numerical results are shown as points marked $+$.
		Results are shown for the first eleven nodes (from $\nodek = 0$ to $\nodek = 10$); results from $\nodek = 1$ are indicated with a thicker blue line and subsequent nodes produce curves further to the right.
		Parameters used are $\lamtotn{} = 1$, $\phin{} = 0.8$, $\musimpn{} = 1.2$, $\nusimpn{} = 0.05$, and $\gamman{} = 0.1$.
		For consistency with the value of $\phin{}$, the initial conditions are $\InfInit = 0.2$ and $\ExpSimpInit = 0.8$.} \label{fig:SEIR-Chain-SusInfAllNodes}
\end{figure}

\subsection{Bounds for SEIR dynamics on a general network}
\label{S:Bounds-SEIR}

We now replicate the argument in \cref{S:Bounds-SIR} to obtain bounds on the solution of generalised SEIR dynamics on a general network.
In this case our starting point is system \eqref{E:GeneralEvol-SEIR} and we assume without loss of generality that $\Rec{\nodek}(0) = 0$ for all nodes.

By analogous arguments to \cref{S:Bounds-SIR}, we observe that 
\begin{equation}
	\InfSus{\nodej}{\nodek} \geq \Sus{\nodek} - \Sus{\nodej} - \uvec \cdot \ExpVec{\nodej} - \RecSus{\nodej}{\nodek}, \label{E:InfSusInequal-SEIR}
\end{equation}
that
\begin{equation}
	\diff{\RecSus{\nodej}{\nodek}}{t} \leq \nuvecn{\nodej} \cdot \ExpVec{\nodej} + \gamman{\nodej} \InfSus{\nodej}{\nodek}, \label{E:RecSusInequal-SEIR}
\end{equation}
and that
\begin{equation}
	-\diff{\Sus{\nodek}}{t} \geq  \lamtotnp{k}{j} \InfSus{\nodej}{\nodek} \label{E:SusEvolInequal-SEIR}
\end{equation}
for any $\nodej \in \Neigh(\nodek)$.

We also replicate some of the analysis from \cref{S:RootedTree-SEIR}. 
We define $\Mn{\nodek}$ as in \eqref{E:MDefn} and we observe that equations \eqref{E:SusEvol-SEIR} and \eqref{E:ExpEvol-SEIR} can be rearranged to give
\begin{equation}
	\begin{bmatrix}
		\diff{\Sus{\nodek}}{t} \\[6pt]
		\diff{\ExpVec{\nodek}}{t}
	\end{bmatrix}
	=
	- \Mn{\nodek}
	\begin{bmatrix}
		\displaystyle \sum_{j \in \Neigh(k)} \lamtotnp{\nodek}{\nodej} \InfSus{\nodej}{\nodek} \\[6pt]
		\ExpVec{\nodek}
	\end{bmatrix}. \label{E:SusAndExpEvolBnd-SEIR}
\end{equation}
Note that \eqref{E:SusAndExpEvolBnd-SEIR} is only valid because $\phivecn{\nodek}$ depends only on $\nodek$ not on the possible sources of infection.
If this were not the case, then it would not be possible to collect the summation terms in the vector on the right hand side of \eqref{E:SusAndExpEvolBnd-SEIR}. 

Repeating the manipulations from \cref{S:RootedTree-SEIR}, we find that 
\begin{equation}
	\nuvecn{\nodej} \cdot \ExpVec{\nodej} = -\nuvecn{\nodej}^T \An{\nodej}^{-1} \phivecn{\nodej} \diff{\Sus{\nodej}}{t} - \nuvecn{\nodej}^T \An{\nodej}^{-1} \diff{\ExpVec{\nodej}}{t}. 	\label{E:NuDotExp-BoundArgument}	
\end{equation}
%and hence
%\begin{equation}
%	\nuvecn{\nodej} \cdot \ExpVec{\nodej} = 
%	- \uvec \cdot \phivecn{\nodej} \diff{\Sus{\nodej}}{t} 
%	+ \mvecn{\nodej} \cdot \phivecn{\nodej} \diff{\Sus{\nodej}}{t}
%	- \uvec \cdot \diff{\ExpVec{\nodej}}{t}
%	+ \mvecn{\nodej} \cdot \diff{\ExpVec{\nodej}}{t}. 
%\end{equation}

Combining \eqref{E:RecSusInequal-SEIR}, \eqref{E:SusEvolInequal-SEIR}, and \eqref{E:NuDotExp-BoundArgument}, we find that 
\begin{equation}
	\diff{\RecSus{\nodej}{\nodek}}{t} \leq  
	 -\nuvecn{\nodej}^T \An{\nodej}^{-1} \phivecn{\nodej} \diff{\Sus{\nodej}}{t} - \nuvecn{\nodej}^T \An{\nodej}^{-1} \diff{\ExpVec{\nodej}}{t}
	- \frac{\gamman{\nodej}}{\lamtotnp{\nodek}{\nodej}} \diff{\Sus{\nodek}}{t}.
\end{equation}
Integrating from $t = 0$ and using the fact that $\RecSus{\nodej}{\nodek}(0) = 0$, we obtain an upper bound on $\RecSus{\nodej}{\nodek}$ that can be substituted into \eqref{E:InfSusInequal-SEIR} and rearranged to obtain
\begin{multline}
	\InfSus{\nodej}{\nodek} \geq 
	\frac{\lamtotnp{\nodek}{\nodej} + \gamman{\nodej}}{\lamtotnp{\nodek}{\nodej}} \Sus{\nodek}
	- (1-\uvec\cdot\phivecn{\nodej})\Sus{\nodej}(t) 
	- \muvecn{\nodej}^T \An{\nodej}^{-1} \left[ \phivecn{\nodej} \Sus{\nodej} + \ExpVec{\nodej} \right]
	\\
	-\nuvecn{\nodej}^T \An{\nodej}^{-1}  \left[ \phivecn{\nodej} \Sus{\nodej}(0) + \ExpVec{\nodej}(0) \right]
	- \frac{\gamman{\nodej}}{\lamtotnp{\nodek}{\nodej}}\Sus{\nodek}(0).
\end{multline}
Since it is also true that $\InfSus{\nodej}{\nodek} \geq 0$, we can use $[x]^+$ as defined in \eqref{E:PosPartDefn-SIR} to obtain a bound on $\InfSus{\nodej}{\nodek}$ analogous to \eqref{E:MainInfSusInequal-SIR}. 
Substituting into \eqref{E:SusEvol-SEIR} then yields 
\begin{multline}
	\diff{\Sus{\nodek}}{t} \leq  
	- \sum_{j \in \Neigh(k)} 
	\Big[(
	- \lamtotnp{\nodek}{\nodej} \nuvecn{\nodej}^T \An{\nodej}^{-1} \left[ \phivecn{\nodej} \Sus{\nodej}(0) + \ExpVec{\nodej}(0) \right]
	- \gamman{\nodej}\Sus{\nodek}(0)
	\\
	+ (\lamtotnp{\nodek}{\nodej} + \gamman{\nodej}) \Sus{\nodek}
	- \lamtotnp{\nodek}{\nodej} (1-\uvec\cdot\phivecn{\nodej}) \Sus{\nodej}(t)
	- \lamtotnp{\nodek}{\nodej} \muvecn{\nodej}^T \An{\nodej}^{-1} \left[ \phivecn{\nodej} \Sus{\nodej} + \ExpVec{\nodej} \right]
	\Big]^{+}. \label{E:MainSusEvol-SEIR-Inequal}
\end{multline}

We note that \eqref{E:MainSusEvol-SEIR-Inequal} depends only on the probabilities of nodes being susceptible or exposed.
Hence, \eqref{E:MainSusEvol-SEIR-Inequal} can be coupled with \eqref{E:ExpEvol-SEIRAlt} to obtain a closed system.
As in \cref{S:Bounds-SIR}, we will use this closed system to show that 
In this case, however, we need to rearrange the system before we can apply the Kamke--M{\"u}ller conditions.

Based on the forms of \eqref{E:MainSusEvol-SEIR-Inequal} and \eqref{E:ExpEvol-SEIRAlt}, we define $\QuuVec{\nodek}(t) = \An{\nodek}^{-1} \left[ \phivecn{\nodek}\Sus{\nodek}(t) + \ExpVec{\nodek}(t)\right]$.
We note that all entries of $\An{\nodek}^{-1}$ are nonnegative and so $\QuuVec{\nodek}$ is nonnegative. 
Rearranging to obtain $\ExpVec{\nodek} = \An{\nodek} \QuuVec{\nodek} - \phivecn{\nodek} \Sus{\nodek}$ and substituting into  \eqref{E:MainSusEvol-SEIR-Inequal} and \eqref{E:ExpEvol-SEIRAlt} then yields
\begin{subequations}
\label{E:SusQuuSystem}
\begin{align}
	\diff{\Sus{\nodek}}{t} &\leq  
	- \sum_{j \in \Neigh(k)} 
	\Big[(
	- \lamtotnp{\nodek}{\nodej} \nuvecn{\nodej} \cdot  \QuuVec{\nodej}(0)
	- \gamman{\nodej}\Sus{\nodek}(0)
	+ (\lamtotnp{\nodek}{\nodej} + \gamman{\nodej}) \Sus{\nodek}
	\notag \\[-12pt]
	& \qquad \qquad \qquad \quad - \lamtotnp{\nodek}{\nodej} (1-\uvec\cdot\phivecn{\nodej}) \Sus{\nodej}(t)
	- \lamtotnp{\nodek}{\nodej} \muvecn{\nodej} \cdot \QuuVec{\nodej}
	\Big]^{+}, \label{E:SusQuuSystem-Sus} \\
	\diff{\QuuVec{\nodek}}{t} &= -\An{\nodek} \QuuVec{\nodek} + \phivecn{\nodek} \Sus{\nodek}. \label{E:SusQuuSystem-Quu}
\end{align}
\end{subequations}

System \eqref{E:SusQuuSystem} is a system of differential inequalities and equations;
as in \cref{S:Bounds-SIR}, we now consider the relationship between the true solutions $\Sus{\nodek}$ and $\QuuVec{\nodek}$ and the approximate solutions $\SusApprox{\nodek}$ and $\QuuVecApprox{\nodek}$ that satisfy the equivalent of \eqref{E:SusQuuSystem} where the inequality in \eqref{E:SusQuuSystem-Sus} is replaced with an equation.
Since the off-diagonal elements of $\An{\nodek}$ are all nonpositive, since $1-\uvec\cdot\phivecn{\nodek} \geq 0$, and since the elements of $\phivecn{\nodek}$ and $\muvecn{\nodek}$ are all nonnegative, this system will satisfy the Kamke--M\"{u}ller conditions and be cooperative.
Hence, we can again apply Lemma 1 from \cite{Simon2018} to conclude that $\SusApprox{\nodek}(t) \geq \Sus{\nodek}(t)$ and that $\QuuVecApprox{\nodek}(t) \geq \QuuVecApprox{\nodek}(t)$ for all $\nodek$ and for all $t$.

While $\QuuVec{\nodek}$ is a useful theoretical construct, we will generally formulate and solve the SEIR rooted-tree approximation using $\ExpVec{\nodek}$ rather than $\QuuVec{\nodek}$.
Using stars to indicate approximate solutions as previously, we use \eqref{E:MainSusEvol-SEIR-Inequal} and \eqref{E:GeneralEvol-SEIRAlt} to obtain the following system as the SEIR rooted-tree approximation:
\begin{subequations}
	\label{E:RootedTreeApprox-SEIR}
	\begin{align}
		\diff{\SusApprox{\nodek}}{t} &=
		- \sum_{j \in \Neigh(k)} 
		\Big[(
		- \lamtotnp{\nodek}{\nodej} \nuvecn{\nodej}^T \An{\nodej}^{-1} \left[ \phivecn{\nodej} \SusApprox{\nodej}(0) + \ExpVecApprox{\nodej}(0) \right]
		- \gamman{\nodej}\SusApprox{\nodek}(0)
		\notag \\[-12pt]
		&\qquad \qquad \qquad 
		+ (\lamtotnp{\nodek}{\nodej} + \gamman{\nodej}) \SusApprox{\nodek}
		- \lamtotnp{\nodek}{\nodej} (1-\uvec\cdot\phivecn{\nodej}) \SusApprox{\nodej}(t)
		\notag \\[-4pt]
		&\qquad \qquad \qquad 
		- \lamtotnp{\nodek}{\nodej} \muvecn{\nodej}^T \An{\nodej}^{-1} \left[ \phivecn{\nodej} \SusApprox{\nodej} + \ExpVecApprox{\nodej} \right]
		\Big]^{+}, \label{E:RootedTreeApprox-SEIR-Sus} \\
		\diff{\ExpVecApprox{\nodek}}{t}  &=  -\phivecn{\nodek} \diff{\SusApprox{\nodek}}{t} - \An{\nodek} \ExpVecApprox{\nodek}, \label{E:RootedTreeApprox-SEIR-Exp} \\
		\diff{\InfApprox{\nodek}}{t} 
		&= 
		-(1 - \uvec\cdot \phivecn{\nodek}) \diff{\SusApprox{\nodek}}{t} - \gamman{\nodek} \InfApprox{\nodek}  + \muvecn{\nodek} \cdot \ExpVecApprox{\nodek}. \label{E:RootedTreeApprox-SEIR-Inf}
	\end{align}
\end{subequations}

Just as \eqref{E:RootedTreeApprox-SIR} is equivalent to \eqref{E:RootedEvol-SIR-Full} for a rooted tree, we can show that \eqref{E:RootedTreeApprox-SEIR} is equivalent to \eqref{E:RootedEvol-SEIR-Full} for a rooted tree.
To see this, we again use the fact that $\SusApprox{\nodej}(t) \geq \SusApprox{\nodek}(t)$ for any $\nodej \in \Neigh(\nodek)$ other than $\nodej = \parn{\nodek}$, and we also use the fact that $\An{\nodek}^{-T} \muvecn{\nodek} + \An{\nodek}^{-T} \nuvecn{\nodek} = \uvec$.
Given that $\SusApprox{\nodek}(t)$ is a decreasing function of $t$, it follows from these observations that 
\begin{equation}
	\nuvecn{\nodej}^T \An{\nodej}^{-1} \phivecn{\nodej} \SusApprox{\nodej}(0) 
	+ (1-\uvec\cdot\phivecn{\nodej}) \SusApprox{\nodej}(t)
	+ \muvecn{\nodej}^T \An{\nodej}^{-1} \phivecn{\nodej} \SusApprox{\nodej}(t)
	\geq \SusApprox{\nodej}(t)
\end{equation}
and hence the term inside the square brackets in \eqref{E:RootedTreeApprox-SEIR-Sus} will be nonpositive whenever $j \neq \parn{\nodek}$.
As a result, \eqref{E:RootedTreeApprox-SEIR} will yield exact solutions for rooted trees without it being necessary to compute the parent of each node.

\begin{figure}
	\centering
	\begin{subfigure}[b]{0.49\textwidth}
		\centering
		\tikzsetnextfilename{SEIR-ER-SusSelectedNodes}
		\begin{tikzpicture}
			\pgfplotstableread{fig-SEIR-ER-SusSelectedNodes-All-A.txt}\mainTable
			\begin{axis}[
				ymin=0,ymax=1,xmin=0,xmax=5,
				xlabel={$t$},
				ylabel={$\Sus{\nodek}$},
				%		xtick={0,0.05,0.1,0.15,0.2},
				%		xticklabels={0,0.05,0.1,0.15,0.2},
				%		xtick={0,0.05,0.1,0.15,0.2,0.25,0.3},
				%		xticklabels={0,0.05,0.1,0.15,0.2,0.25,0.3},
				%	ytick={0},
				%	extra x ticks={0.5,1},
				%	extra y ticks={0.5,1},
				%	extra x tick labels={$\tfrac{1}{2}$, $1$}, 
				%	extra y tick labels={$\tfrac{1}{2}$, $1$}, 
				width={\textwidth},
				height={0.8\textwidth}
				]
				\addplot[color=black,mark=o] table[x index = 0, y index = 1] from \mainTable;
				\addplot[color=black,mark=o,dashed,mark options=solid] table[x index = 0, y index = 5] from \mainTable;
				\addplot[color=blue,mark=+] table[x index = 0, y index = 2] from \mainTable;
				\addplot[color=blue,mark=+,dashed,mark options=solid] table[x index = 0, y index = 6] from \mainTable;
				\addplot[color=red,mark=square] table[x index = 0, y index = 3] from \mainTable;
				\addplot[color=red,mark=square,dashed,mark options=solid] table[x index = 0, y index = 7] from \mainTable;
				\addplot[color=olive,mark=triangle] table[x index = 0, y index = 4] from \mainTable;
				\addplot[color=olive,mark=triangle,dashed,mark options=solid] table[x index = 0, y index = 8] from \mainTable;
				%	\foreach \figLineMaker in {2,3,...,6} {
				%		\addplot[color=black] table[x index = 0, y index = \figLineMaker] from \mainTable;
				%	}
				%	\foreach \figLineMaker in {7,8,...,12} {
				%		\addplot[color=black, dashed] table[x index = 0, y index = \figLineMaker] from \mainTable;
				%	}
			\end{axis}
		\end{tikzpicture}
		%	\caption{
		%		Probabilities of selected nodes (indicated using different colours and mark styles) in an Erd\H{o}s--R\'{e}nyi graph being susceptible at time $t$, according to the rooted-tree approximation (continuous lines) and the average of $10^5$ Gillespie simulations (dashed lines).} 
		\caption{}
		\label{fig:SEIR-ER-SusSelectedNodes}
	\end{subfigure}
	\begin{subfigure}[b]{0.49\textwidth}
		\centering
		\tikzsetnextfilename{SEIR-ER-AvgSusInf}
		\begin{tikzpicture}
			\pgfplotstableread{fig-SEIR-ER-AvgSusExpInfRec-All-A.txt}\mainTable
			\begin{axis}[
				ymin=0,ymax=1,xmin=0,xmax=12,
				xlabel={$t$},
				ylabel={Averaged $\Sus{\nodek}$, $\ExpSimp{\nodek}$, $\Inf{\nodek}$, $\Rec{\nodek}$},
				%			xtick={0,0.05,0.1,0.15,0.2,0.25,0.3},
				%			xticklabels={0,0.05,0.1,0.15,0.2,0.25,0.3},
				%	ytick={0},
				%	extra x ticks={0.5,1},
				%	extra y ticks={0.5,1},
				%	extra x tick labels={$\tfrac{1}{2}$, $1$}, 
				%	extra y tick labels={$\tfrac{1}{2}$, $1$}, 
				width={\textwidth},
				height={0.8\textwidth}
				]
				\addplot[color=black,ultra thick] table[x index = 0, y index = 1] from \mainTable;
				\addplot[color=black,ultra thick,dashed] table[x index = 0, y index = 5] from \mainTable;
				\addplot[color=olive,thick] table[x index = 0, y index = 2] from \mainTable;
				\addplot[color=olive,dashed,thick] table[x index = 0, y index = 6] from \mainTable;
				\addplot[color=blue] table[x index = 0, y index = 3] from \mainTable;
				\addplot[color=blue,dashed] table[x index = 0, y index = 7] from \mainTable;
				\addplot[color=red,thin] table[x index = 0, y index = 4] from \mainTable;
				\addplot[color=red,dashed,thin] table[x index = 0, y index = 8] from \mainTable;
				%	\foreach \figLineMaker in {2,3,...,6} {
				%		\addplot[color=black] table[x index = 0, y index = \figLineMaker] from \mainTable;
				%	}
				%	\foreach \figLineMaker in {7,8,...,12} {
				%		\addplot[color=black, dashed] table[x index = 0, y index = \figLineMaker] from \mainTable;
				%	}
			\end{axis}
		\end{tikzpicture}
		%	\caption{Averaged probabilities over all nodes of nodes being susceptible (black/thick lines), infectious (blue/medium lines), or recovered (red/thin lines) in an Erd\H{o}s--R\'{e}nyi graph according to the rooted-tree approximation (continuous lines) and the average of $10^5$ Gillespie simulations (dashed lines).}
		\caption{}
		\label{fig:SEIR-ER-AvgSusExpInf}
	\end{subfigure}

	\caption{\small 
		Comparisions of the rooted-tree approximation in \eqref{E:RootedTreeApprox-SEIR} with simulation results from the average of $10^5$ Gillespie algorithm simulations of the full stochastic SEIR model for an Erd{\H o}s--R\'enyi random graph with 100 nodes and probability of connection 0.05.
		Subfigure (a) shows $\Sus{\nodek}$ for four different nodes: results from the rooted-tree approximation are shown as continuous lines and results from Gillespie simulations are shown as dashed lines;
		different nodes are distinguished using different colours and marker styles.
		Subfigure (b) show $\Sus{\nodek}$ (very thick black lines), $\ExpSimp{\nodek}$ (thick olive lines), $\Inf{\nodek}$ (medium thickness blue lines) and $\Rec{\nodek}$ (thin red lines) averaged over all nodes in the network:
		results from the rooted-tree approximation are shown as continuous lines and results from Gillespie simulations are shown as dashed lines.
		Parameters used are $\lamtotn{} = 1$, $\phin{} = 0.8$, $\musimpn{} = 1.2$, $\nusimpn{} = 0.05$, and $\gamman{} = 0.1$.
		There is a single node that is infectious at $t=0$ and all other nodes are susceptible.
	} \label{fig:SEIR-Comparisons} %\label{fig:SIR-AlmostTree}
\end{figure}
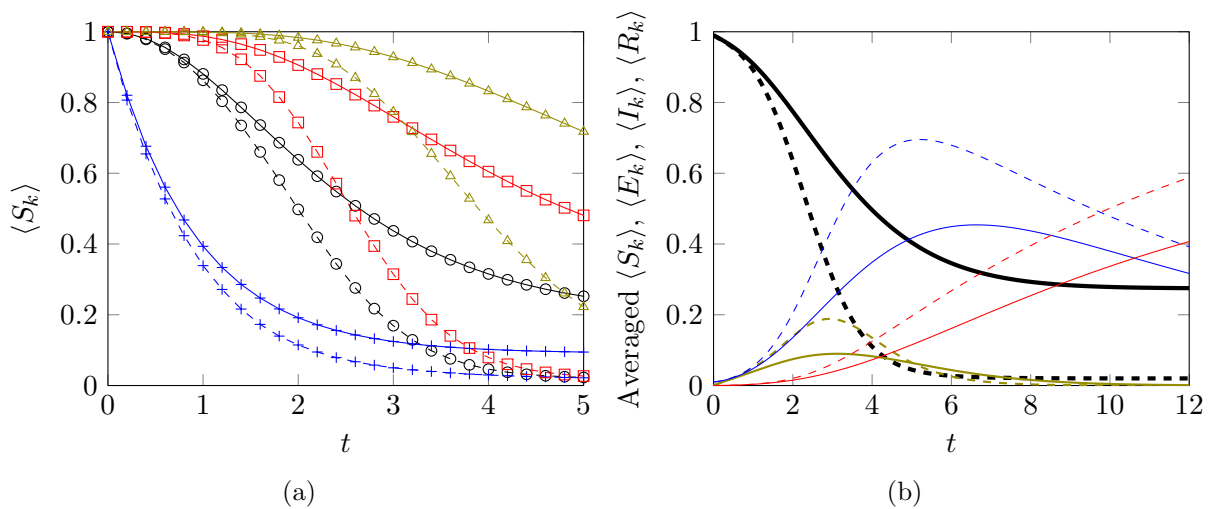

\cref{fig:SEIR-Comparisons} is analogous to \cref{fig:SIR-Comparisons} and it enables equivalent conclusions to be drawn.
\cref{fig:SEIR-Comparisons} shows comparisons of the rooted-tree approximation \eqref{E:RootedTreeApprox-SEIR} with estimates of the true solution obtained from averaging $10^5$ simulations using the Gillespie algorithm for an ER random graph.
Code is available at \url{https://github.com/cameronlhall/rootedtreeapprox} and the ER graph used to generate \cref{fig:SEIR-Comparisons} is different from the ER graph used in \cref{fig:SIR-Comparisons}.

From \cref{fig:SEIR-ER-SusSelectedNodes} we verify that the rooted-tree approximation gives an upper bound on $\Sus{\nodek}$ for the nodes $\nodek$ illustrated.
From \cref{fig:SEIR-ER-AvgSusExpInf}, we see that there is a reasonably large difference between the true solution (dashed lines) and the rooted-tree approximation (continuous lines) and so once again the bounds provided by \eqref{E:RootedTreeApprox-SEIR} are not generally very tight.

\section{Discussion and conclusions}
\label{S:Discussion}

In this paper, we have developed and analysed a new approximation method, the rooted-tree approximation, that can be applied to SIR and generalised SEIR models on networks.
In the case of a tree with a unique initially-infected node, our approximation is exact and leads to a partially-decoupled system of linear differential equations for the node-state probabilities.
As demonstrated in \cref{S:Solutions-SIR}, we can obtain explicit closed-form solutions for the node state probabilities for SIR models and, in theory, equivalent results can also be obtained for SEIR models.

Since the pair-based and message-passing approximations are both exact on \textit{all} trees (not just rooted trees) but closed-form solutions for these are not well known, it is instructive to compare our system \eqref{E:RootedEvol-SIR-Full} with appropriate rooted tree simplifications of the pair-based SIR approximation in \cite{Sharkey2015} and the message-passing SIR approximation in \cite{Karrer2010}.
For the pair-based approximation (\eg, system (3) in \cite{Sharkey2015}), we find that we can use proof by induction from the leaves to the root to show that $\InfSus{\nodek}{\parn{\nodek}} = 0$ on a rooted tree. 
Subsequently, we can use the fact that  $\SusSus{\parn{\nodek}}{\nodek} = \Sus{\nodek}$ to convert the remaining equations of the pair-based approximation into a linear system equivalent to \eqref{E:RootedEvol-SIR-Full}.

For the message-passing model in \cite{Karrer2010} applied to a rooted tree, we can work from the leaves to the root to show that $H^{\parn{\nodek}\leftarrow\nodek} = 0$ and then work back out from the root to the leaves to obtain expressions for $H^{\nodek\leftarrow\parn{\nodek}}$ that are analogous to an integrated form of our system \eqref{E:RootedEvol-SIR-Full}. 
As a result, we find that the explicit solutions in \eqref{E:ClosedFormSolution} could have been obtained from the pair-based or message-passing approximations; 
while we believe that this is the first time that these explicit solutions have been reported, they are consistent with---and theoretically obtainable from---established results in the existing literature.

One important feature of our rooted-tree approximation is that it provides upper bounds on $\Sus{\nodek}$ at every node.
This is an important strength of our method since it provides a contrast from other methods that yield lower bounds on $\Sus{\nodek}$.
One promising avenue for further research is to combine the rooted-tree approximation with other approximations in order to obtain better estimates of node-state probabilities.
Such hybrid approximations are likely to be more practical than the rooted-tree approximation because the bounds on $\Sus{\nodek}$ are rarely very tight.
As we see from \cref{fig:SIR-Comparisons,fig:SEIR-Comparisons}, there are often large differences between the node-state probabilities obtained from the rooted-tree approximation and estimates of the true node-state probabilities based on Gillespie algorithm simulations.

Another limitation of the rooted-tree approximation is that it is reliant on assumptions that there can be no return to a susceptible state and that there can only be one variety of infectious state.
Both of these assumptions are necessary in order to express $\InfSus{\parn{\nodek}}{\nodek}$, and hence the rate of infection, in terms of a linear combination of the node-state probabilities and $\RecSus{\parn{\nodek}}{\nodek}$ for rooted trees.
One avenue for further research would be to explore whether the rooted-tree approximation can be extended to SIRS and SEIRS models or SIR models with multiple infecious states.
Perhaps this would involve developing new approximations that are not exact on rooted trees but would still provide a consistent upper bound on $\Sus{\nodek}$, analogous to the $W(x,y) = \min(x,y)$ approximation for SIS models introduced in \cite{Simon2018}.

Overall, the rooted-tree approximation presented in this paper is a new way of analysing SIR and SEIR dynamics on networks that has advantages and disadvantages over existing methods.
The principal strengths of the rooted-tree approximation are that it is simple (leading to a cooperative, piecewise-linear system of equations for node-state probabilities), that it yields exact closed-form solutions in certain situations, and that it yields upper bounds on $\Sus{\nodek}$ in contrast with the lower bounds provided by other approximations. 
The principal weakness of the rooted-tree approximation is that the bounds on $\Sus{\nodek}$ are not very tight unless the underlying network is a tree with a single initially-infected node.
Despite this limitation, the simplicity of the rooted-tree approximation means that it has the potential to be a useful tool in developing new computational methods for analysing contagion dynamics on networks.

\bibliographystyle{siamplain}
\bibliography{contagionBib}

%\appendix
%\section{Explicit solution of SIR on chain}
%\label{S:SolutionMethod}

\end{document}